\documentclass[a4paper]{article}
\usepackage[usenames,svgnames,dvipsnames,table]{xcolor}
\usepackage{graphicx}
\usepackage{morefloats}
\usepackage{pdflscape}
\usepackage{pdfpages}
\usepackage{booktabs}
\usepackage{comment}
\usepackage{textcomp}
\usepackage[utf8]{inputenc}
\usepackage{natbib}
\setcitestyle{authoryear}
\usepackage{enumitem}
\usepackage{url}
\usepackage{balance} 
\usepackage{amsfonts}
\usepackage{multirow}
\usepackage[colorlinks]{hyperref}
\usepackage{amsmath}
\usepackage[framemethod=tikz]{mdframed}
\newcommand{\alg}[1]{\texttt{#1}}
\newcommand{\myvec}[1]{\mathbf{#1}}
\newcommand{\ie}{i.e.}

\definecolor{NavyBlue}{RGB}{35,35,142}
\definecolor{RawSienna}{RGB}{199,97,20}
\hypersetup{
	colorlinks,%
	citecolor=NavyBlue,%
	filecolor=NavyBlue,%
	linkcolor=RawSienna,%
	urlcolor=NavyBlue
}

\usepackage[linesnumbered,ruled,vlined]{algorithm2e}
\usepackage{algorithmicx}
\usepackage{algpseudocode}

\newcommand{\hv}{$I^{-}_{H}$}
\newcommand{\gd}{$I_{GD}$}
\newcommand{\igd}{$I_{IGD}$}
\newcommand{\epsind}{$I^{1}_{\epsilon+}$}

\SetKwFunction{ND}{ND}

\title{\Large \textbf{BMOBench}: Black-Box Multi-Objective Optimization Benchmarking Platform\\
	\vspace{1em}
\large {\it
	Supplementary Materials for
Multi-Objective Simultaneous Optimistic
Optimization}
}

\begin{document}
\author{Abdullah Al-Dujaili, Suresh Sundaram}
	\maketitle

	
This document briefly describes the Black-Box Multi-Objective Optimization Benchmarking (\textbf{BMOBench}) platform. It presents the test problems,  evaluation procedure, and experimental setup. To this end, the BMOBench is demonstrated by comparing recent multi-objective solvers from the literature, namely  \alg{SMS-EMOA}~\citep{beume2007sms}, \alg{DMS}~\citep{custodio2011direct}, and \alg{MO-SOO}~\citep{ash-mosoo-15}.

\section{Test Problems}

\begin{table}[tb]
	\renewcommand{\arraystretch}{1.2}
	\begin{center}
		\resizebox{\textwidth}{!}{
	\begin{tabular}{
			>{\centering\arraybackslash}p{0.25in}p{2.5in}>{\centering\arraybackslash}p{0.15in}>{\centering\arraybackslash}p{0.25in}
			>{\centering\arraybackslash}p{0.25in}>{\centering\arraybackslash}p{0.25in}>{\centering\arraybackslash}p{0.25in}|
			>{\centering\arraybackslash}p{0.25in}p{2.5in}>{\centering\arraybackslash}p{0.25in}>{\centering\arraybackslash}p{0.25in}
			>{\centering\arraybackslash}p{0.25in}>{\centering\arraybackslash}p{0.25in}>{\centering\arraybackslash}p{0.25in}|
			>{\centering\arraybackslash}p{0.25in}p{2.5in}>{\centering\arraybackslash}p{0.25in}>{\centering\arraybackslash}p{0.25in}
			>{\centering\arraybackslash}p{0.25in}>{\centering\arraybackslash}p{0.25in}>{\centering\arraybackslash}p{0.25in}} 
		\toprule
		\# & Problem Name & n & m & \textbf{D} & \textbf{S} & \textbf{M} & \# & Problem Name & n & m & \textbf{D} & \textbf{S} & \textbf{M} & \# & Problem Name & n & m & \textbf{D} & \textbf{S} & \textbf{M} \\
		\toprule
		1  & BK1 \citep{25_p} 							& 2  & 2 &  L & S & U   				&	35 & I5 \citep{24_p} 				& 8  & 3 &  H & NS & U				&	68 & MOP3 \citep{25_p} 		& 2  & 2 & L & $\times$& $\times$ \\
		2  & CL1 \citep{6_p}                			& 4  & 2 &  L & $\times$ & $\times$		&	36 & IKK1 \citep{25_p} 				& 2  & 3 &  L & $\times$&U			&	69 & MOP4 \citep{25_p} 		& 3  & 2 & L & S&$\times$\\
		3  & Deb41 \citep{15_p}							& 2  & 2 &  L & $\times$ & $\times$		&	37 & IM1 \citep{25_p} 				& 2  & 2 &  L & $\times$&U			&	70 & MOP5 \citep{25_p} 		& 2  & 3 & L & NS&$\times$\\
		4  & Deb512a \citep{15_p} 						& 2  & 2 &  L & $\times$ & $\times$		&	38 & Jin1 \citep{48_p} 				& 2  & 2 &  L &$\times$ &U			&	71 & MOP6 \citep{25_p} 		& 2  & 2 & L & S&$\times$\\
		5  & Deb512b \citep{15_p} 						& 2  & 2 &  L & $\times$&$\times$		&	39 & Jin2 \citep{48_p} 				& 2  & 2 &  L &$\times$ &U			&	72 & MOP7 \citep{25_p} 		& 2  & 3 & L & $\times$&U\\
		6  & Deb512c \citep{15_p} 						& 2  & 2 &  L & $\times$&$\times$		&	40 & Jin3 \citep{48_p} 				& 2  & 2 &  L &$\times$ &U			&	73 & OKA1 \citep{39_p} 		& 2  & 2 & L & $\times$&$\times$\\
		7  & Deb513 \citep{15_p} 						& 2  & 2 &  L & $\times$&$\times$		&	41 & Jin4 \citep{48_p} 				& 2  & 2 &  L &$\times$ &U			&	74 & OKA2 \citep{39_p} 		& 3  & 2 & L & $\times$&$\times$\\
		8  & Deb521a \citep{15_p}					 	& 2  & 2 &  L & $\times$&$\times$		&	42 & Kursawe \citep{31_p} 			& 3  & 2 &  L &$\times$ &$\times$	&	75 & QV1 \citep{25_p} 		& 10 & 2 & H & S&M\\
		9  & Deb521b \citep{15_p} 						& 2  & 2 &  L & $\times$&$\times$		&	43 & L1ZDT4 \citep{18_p} 			& 10 & 2 &  H &$\times$ &$\times$	&	76 & Sch1 \citep{25_p} 		& 1  & 2 & L & $\times$&$\times$\\
		10 & Deb53 \citep{15_p} 						& 2  & 2 &  L & $\times$&$\times$		&	44 & L2ZDT1 \citep{18_p} 			& 30 & 2 &  H &$\times$ &$\times$ 	&	77 & SK1 \citep{25_p} 		& 1  & 2 & L & S&M\\
		11 & DG01 \citep{25_p} 							& 1  & 2 &  L & $\times$ & M			&	45 & L2ZDT2 \citep{18_p} 			& 30 & 2 &  H &$\times$ &$\times$	&	78 & SK2 \citep{25_p} 		& 4  & 2 & L & $\times$&$\times$\\
		12 & DPAM1 \citep{25_p} 						& 10 & 2 &  H & NS & $\times$			&	46 & L2ZDT3 \citep{18_p} 			& 30 & 2 &  H &$\times$ &$\times$	&	79 & SP1 \citep{25_p} 		& 2  & 2 & L & NS&U\\
		13 & DTLZ1 \citep{19_p} 						& 7  & 3 &  H & $\times$&M				&	47 & L2ZDT4 \citep{18_p} 			& 30 & 2 &  H &$\times$ &$\times$ 	&	80 & SSFYY1 \citep{25_p}    & 2  & 2 & L & S&U\\
		14 & DTLZ1n2 \citep{19_p} 						& 2  & 2 &  L & $\times$&M				&	48 & L2ZDT6 \citep{18_p} 			& 10 & 2 &  H &$\times$ &$\times$	&	81 & SSFYY2 \citep{25_p}    & 1  & 2 & L & $\times$&$\times$\\
		15 & DTLZ2 \citep{19_p} 						& 12 & 3 &  H & $\times$&U				&	49 & L3ZDT1 \citep{18_p} 			& 30 & 2 &  H &$\times$ &$\times$	&	82 & TKLY1 \citep{25_p} 	& 4  & 2 & L & $\times$&$\times$\\
		16 & DTLZ2n2 \citep{19_p} 						& 2  & 2 &  L & $\times$&U				&	50 & L3ZDT2 \citep{18_p} 			& 30 & 2 &  H &$\times$ &$\times$	&	83 & VFM1 \citep{25_p} 		& 2  & 3 & L & S&U\\
		17 & DTLZ3 \citep{19_p} 						& 12 & 3 &  H & $\times$&M				&	51 & L3ZDT3 \citep{18_p} 			& 30 & 2 &  H &$\times$ &$\times$ 	&	84 & VU1 \citep{25_p} 		& 2  & 2 & L & S&U\\
		18 & DTLZ3n2 \citep{19_p} 						& 2  & 2 &  L & $\times$&M				&	52 & L3ZDT4 \citep{18_p} 			& 30 & 2 &  H &$\times$ &$\times$	&	85 & VU2 \citep{25_p} 		& 2  & 2 & L & S&U\\
		19 & DTLZ4 \citep{19_p} 						& 12 & 3 &  H & $\times$&U				&	53 & L3ZDT6 \citep{18_p} 			& 10 & 2 &  H &$\times$ &$\times$	&	86 & WFG1 \citep{25_p} 		& 8  & 3 & H & S&U\\
		20 & DTLZ4n2 \citep{19_p} 						& 2  & 2 &  L & $\times$&U				&	54 & LE1 \citep{25_p} 				& 2  & 2 &  L &S &U					&	87 & WFG2 \citep{25_p} 		& 8  & 3 & H & NS&$\times$\\
		21 & DTLZ5 \citep{19_p} 						& 12 & 3 &  H & $\times$&U				&	55 & lovison1 \citep{33_p} 			& 2  & 2 &  L &$\times$ &$\times$	&	88 & WFG3 \citep{25_p} 		& 8  & 3 & H & NS&U\\
		22 & DTLZ5n2 \citep{19_p} 						& 2  & 2 &  L & $\times$&U				&	56 & lovison2 \citep{33_p} 			& 2  & 2 &  L &$\times$ &$\times$	&	89 & WFG4 \citep{25_p} 		& 8  & 3 & H & S&M\\
		23 & DTLZ6 \citep{19_p} 						& 22 & 3 &  H & $\times$&U				&	57 & lovison3 \citep{33_p} 			& 2  & 2 &  L &$\times$ &$\times$	&	90 & WFG5 \citep{25_p} 		& 8  & 3 & H & S&$\times$\\
		24 & DTLZ6n2 \citep{19_p} 						& 2  & 2 &  L & $\times$&U				&	58 & lovison4 \citep{33_p} 			& 2  & 2 &  L &$\times$ &$\times$	&	91 & WFG6 \citep{25_p} 		& 8  & 3 & H & NS&U\\
		25 & ex005 \citep{26_p} 						& 2  & 2 &  L & $\times$&U				&	59 & lovison5 \citep{33_p} 			& 3  & 3 &  L &$\times$ &$\times$	&	92 & WFG7 \citep{25_p} 		& 8  & 3 & H & S&U\\
		26 & Far1 \citep{25_p} 							& 2  & 2 &  L & NS&	M					&	60 & lovison6 \citep{33_p} 			& 3  & 3 &  L &$\times$ &$\times$	&	93 & WFG8 \citep{25_p} 		& 8  & 3 & H & NS&U\\
		27 & FES1 \citep{25_p} 							& 10 & 2 &  H & S&	U					&	61 & LRS1 \citep{25_p} 				& 2  & 2 &  L &S &U					&	94 & WFG9 \citep{25_p} 		& 8  & 3 & H & NS&$\times$\\
		28 & FES2 \citep{25_p} 							& 10 & 3 &  H & S&	U					&	62 & MHHM1 \citep{25_p} 			& 1  & 3 &  L &$\times$ &U			&	95 & ZDT1 \citep{49_p} 		& 30 & 2 & H & S&U\\
		29 & FES3 \citep{25_p} 							& 10 & 4 &  H & S&	U					&	63 & MHHM2 \citep{25_p} 			& 2  & 3 &  L &S &		U			&	96 & ZDT2 \citep{49_p} 		& 30 & 2 & H & S&U\\
		30 & Fonseca \citep{21_p} 						& 2  & 2 &  L & S& U				&	64 & MLF1 \citep{25_p} 				& 1  & 2 &  L &$\times$ & M			&	97 & ZDT3 \citep{49_p} 		& 30 & 2 & H & S&$\times$\\
		31 & I1 \citep{24_p} 							& 8  & 3 &  H &S &	U					&	65 & MLF2 \citep{25_p} 				& 2  & 2 &  L &NS &M				&	98 & ZDT4 \citep{49_p} 		& 10 & 2 & H & S&$\times$\\
		32 & I2 \citep{24_p} 							& 8  & 3 &  H &NS &	U					&	66 & MOP1 \citep{25_p} 				& 1  & 2 &  L &S &U					&	99 & ZDT6 \citep{49_p} 		& 10 & 2 & H & S&M\\
		33 & I3 \citep{24_p} 							& 8  & 3 &  H &NS &	U					&	67 & MOP2 \citep{25_p} 				& 4  & 2 &  L &S &U					&  100 & ZLT1 \citep{25_p} 		& 10 & 3 & H & S& U\\
		34 & I4 \citep{24_p} 							& 8  & 3 &  H &NS &	 U   &      & 				   				&    &	 &    &	&   &		& 					    &    &   & & & \\
		\bottomrule	
	\end{tabular}
}
	\end{center}
	\caption{Test problems definition and properties. Symbols: \textbf{D} : dimensionality $\in\{$L : low-dimensionality, H : high-dimensionality$\}$; \textbf{S} : separability $\in\{$S : separable, NS : non-separable$\}$; \textbf{M} : modality $\in\{$U : uni-modal, M : multi-modal$\}$; $\times$ : uncategorized/mixed.}\label{tbl:problem_def}
\end{table}

One-hundred multi-objective optimization problems from the literature are selected.\footnote{
	retrieved from \url{ http://www.mat.uc.pt/dms}.} These problems have simple bound constraints, that is to say, $\mathcal{X}=[\myvec{l},\myvec{u}]\subset \mathbb{R}^n$, where $\myvec{u}\succeq\myvec{l}$. Table~\ref{tbl:problem_def} presents a brief list of these problems with number of dimensions/objectives. In order to have a better understanding of the algorithm strength/weakness, the benchmark problems are categorized (wherever possible) according to three key characteristics, namely \textit{dimensionality}: low- or high-dimension decision space, \textit{separability}: separable or non-separable objectives, and \textit{modality}: uni-modal or multi-modal objectives. Each of these attributes imposes a different challenge in solving an \textsc{MOO} problem~\citep{25_p}.

\section{Evaluation Budget}
\label{sec:eval_budget}
\alg{MO-SOO} is a deterministic algorithm producing the same approximation set in each run of the algorithm for a given problem, whereas the approximation sets produced by the compared stochastic algorithms: \alg{DMS} and \alg{SMS-EMOA} can be different every time they are run for a given problem. In practice, stochastic algorithms are run several times per problem. To this end and to ensure a fair comparison, given a computational budget of $v$ function evaluations per run, the stochastic algorithms are allocated 10 runs per problem instance. On the other hand, the deterministic algorithms are run once per problem instance with the accumulated $10\times v$ function evaluations.

In our experiments,  the evaluation budget $v$ is made proportional to the search space dimension~$n$ and is set to $10^2\cdot n$. The overall computational budget used by an algorithm on BMOBench is the product of the evaluation budget per run, the number of problems, and the number of runs per problem.

With $n=2$, for instance, the overall computational budget used by \alg{MO-SOO} on BMOBench is $10^3\cdot 2 \cdot 100\cdot 1= 2\times 10^5$ function evaluations. Each of the other algorithms uses also a computational budget of $10^2\cdot 2 \cdot 100\cdot 10= 2\times 10^5$ function evaluations.
	\section{Benchmark Procedure}
	\label{sec:perf proc}

			\begin{table}[tb]
				\begin{center}
					\renewcommand{\arraystretch}{1.5}
					\resizebox{0.9\textwidth}{!}{
						\begin{small}
							\begin{tabular}{@{}l@{}>{\centering\arraybackslash}p{1.8in}@{}>{\centering\arraybackslash}p{1.8in}@{}>{\centering\arraybackslash}p{1.0in}@{}}
								
								\toprule
								Quality Indicator ($I$) & Pareto-Compliant	& Reference Set Required		& Target  \\ \toprule
								Hypervolume Difference (\hv) &	Yes&	Yes	 &	Minimize\\
								Generational Distance (\gd) &	No&	Yes&		Minimize\\
								Inverted Generational Distance (\igd ) &	No&	Yes&		Minimize\\
								Additive $\epsilon$-Indicator (\epsind) &	Yes&	Yes&		Minimize \\
								\bottomrule		
							\end{tabular}
						\end{small}}
					\end{center}
					\caption{Employed Quality Indicators. Adapted from \citep{hadka2012moea} \cite[for more details, see][]{on_metric_performance,coello2002evolutionary}. }\label{tbl:metrics}
				\end{table}

	Similar to~\citep{brockhoff:hal-01146741}, a set of targets are defined for  each problem in terms of four popular quality indicators~\citep{knowles_tutorial_indicator,zitzler2003performance} listed in Table~\ref{tbl:metrics}. 
	A solver (algorithm) is then evaluated based on its runtime with respect to each target: the number of function evaluations used until the target is reached. We present the recorded runtime values in terms of \emph{data profiles} \citep{more2009benchmarking}. A data profile can be regarded as an empirical cumulative
	distribution function of the observed number of function evaluations in which the y-axis tells how
	many targets---over the set of problems and quality indicators---have been reached by each algorithm for a
	given evaluation budget (on the x-axis). Mathematically, a data profile for a solver $s$ on a problem class $P$ has the form
	\begin{equation}
	d_s(\alpha) = \frac{1}{|P|} \bigg|\bigg\{p \in P\;\Big|\; \frac{t_{p,s}}{n_p} \leq \alpha \bigg\}\bigg|\;,\nonumber
	\end{equation}
	where $t_{p,s}$ is the observed runtime of solver $s$ on solving problem $p$ (hitting a target) over a decision space $\mathcal{X}\subseteq \mathbb{R}^{n_p}$. The data profile approach captures several benchmarking aspects, namely the convergence behavior over time rather than a fixed budget, which can as well be aggregated over problems of similar category~\citep[see, for more details,][]{brockhoff:hal-01146741}. In our experiments, 70 linearly spaced values in the logarithmic scale from $10^{-0.8}$ to $10^{-3}$ and from $10^{-0.1}$ to $10^{-2}$ were used as targets for~\hv,~\gd, and~\igd; and~\epsind, respectively.
	
	The~\hv,~\gd,~\igd~and~\epsind values are computed for each algorithm at any point of its run based on the set of all (normalized) non-dominated vectors found so far---\ie, the \emph{archive}---with respect to a (normalized) \emph{reference} set~$R\in\Omega$. 
	  We computed the reference set for calculating the quality indicators by aggregating(union) the approximation sets generated by the evolutionary algorithms used in~\citep{custodio2011direct}.
	
	As mentioned in the previous section, we aim to provide a fair comparison between deterministic and stochastic solvers and accommodate the multiple-run practice for stochastic algorithms, at the same time. This has been reflected in the evaluation budget allocation (see Section~\ref{sec:eval_budget}). Likewise, we need to adapt the data profiles. To this end, given a problem instance and for each one of the stochastic solvers, we consider the \emph{best} reported runtime for each target from the solver's 10 runs, rather than the mean value. With this setting in hand, the data profile of \alg{MO-SOO} at $10^3$ function evaluations, for instance, can be compared to that of \alg{SMS-EMOA} at $10^2$ function evaluations.

\section{Results}

\begin{figure}[h!]
	\centering
	\includegraphics[width=0.45\textwidth]{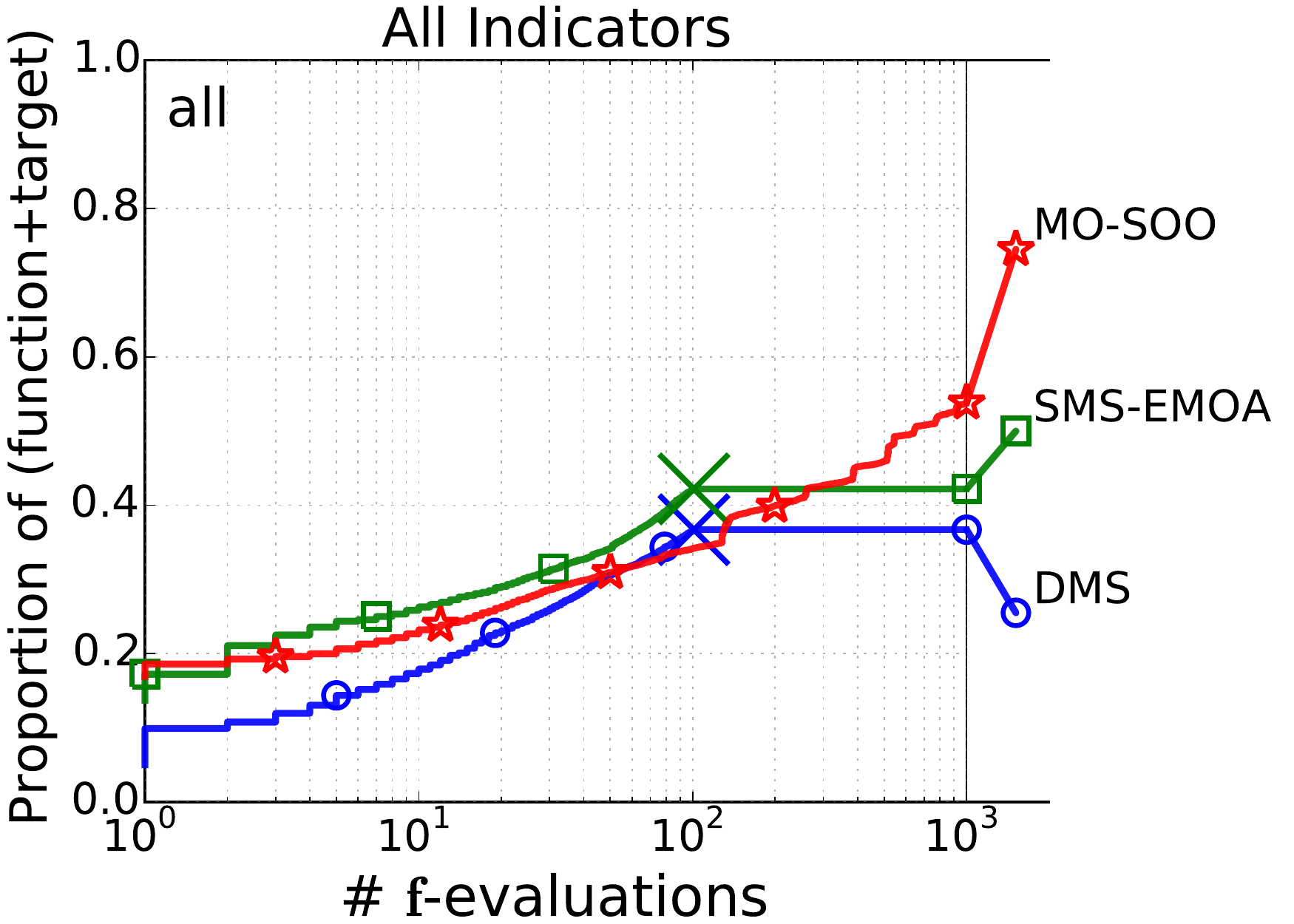}
	\caption{Data profiles aggregated over all the problems across all
		the quality indicators computed for each of the compared algorithms. The symbol × indicates the maximum number of function
		evaluations.
	}
	\label{fig:all_ind}
\end{figure}

\begin{figure*}[h!]
	\begin{center}
		\renewcommand{\arraystretch}{1.1}
		\resizebox{0.95\textwidth}{!}{
			\begin{tabular}{cccc}
				\toprule
				\multicolumn{4}{c}{\textbf{Hypervolume} (HV)}\\
				\includegraphics[width=7cm, trim = 0mm 0mm 0mm 0mm, clip]{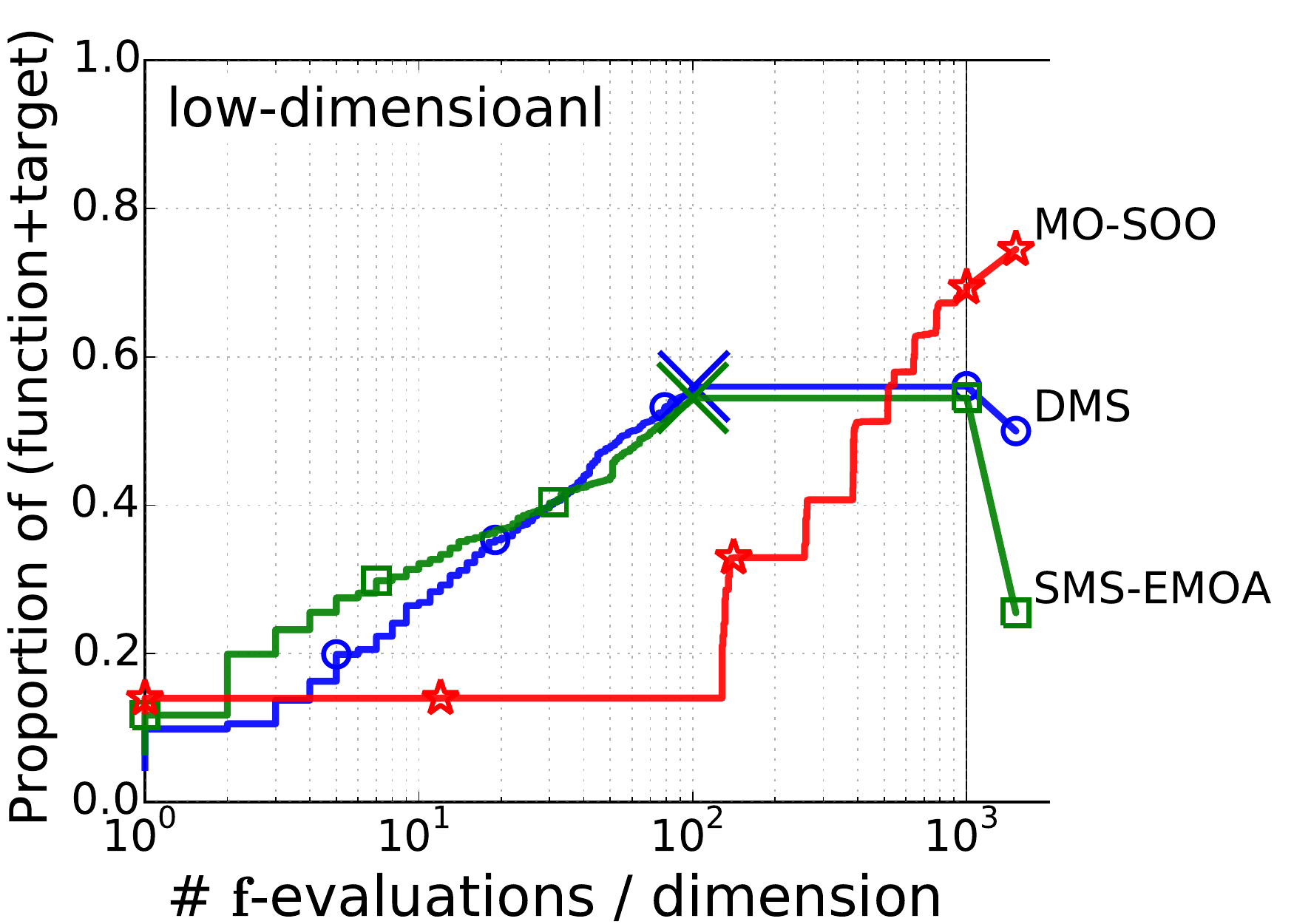}&
				\includegraphics[width=7cm, trim = 0mm 0mm 0mm 0mm, clip]{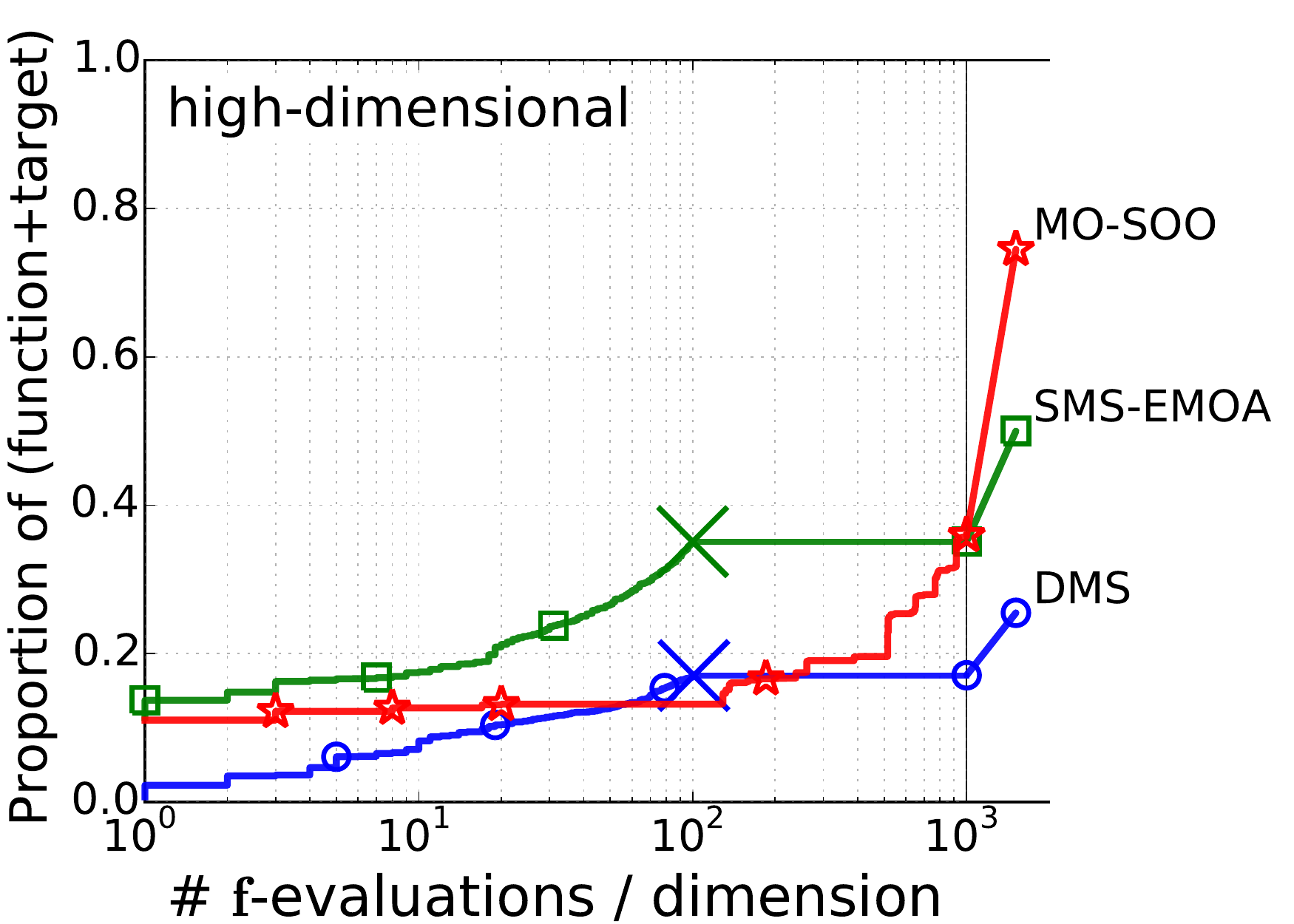}&
				\includegraphics[width=7cm, trim = 0mm 0mm 0mm 0mm, clip]{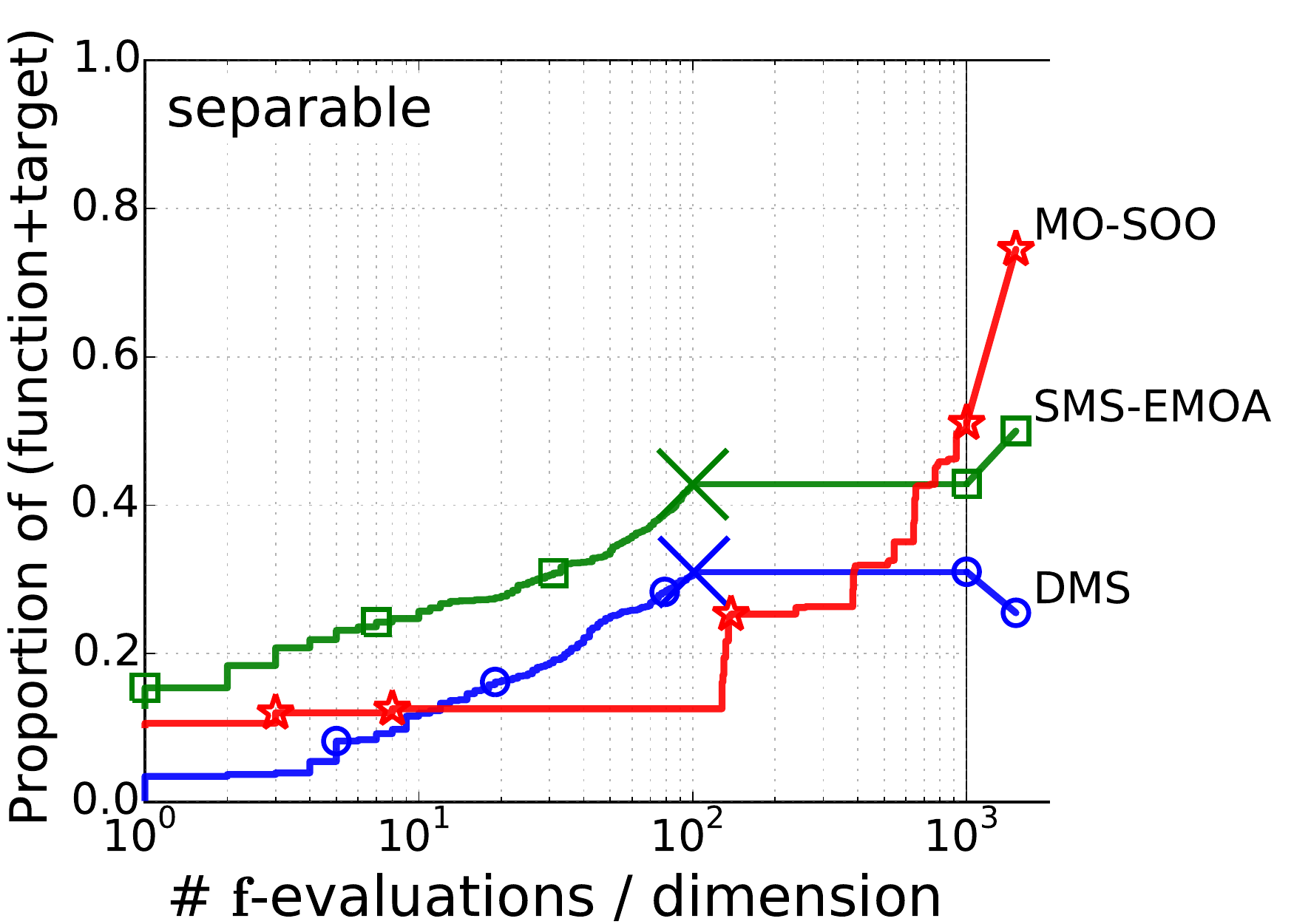}&
				\includegraphics[width=7cm, trim = 0mm 0mm 0mm 0mm, clip]{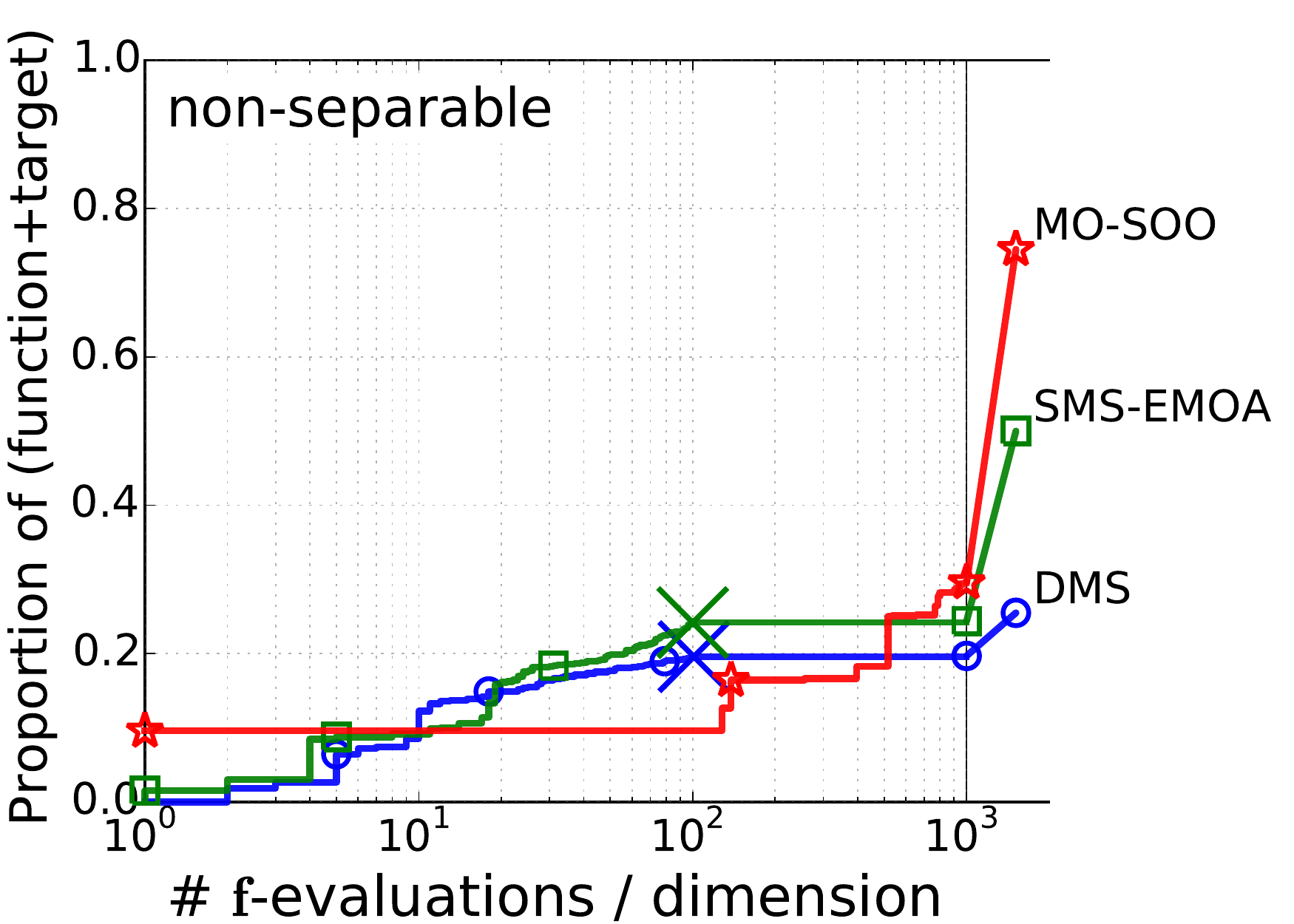}\\
				\includegraphics[width=7cm, trim = 0mm 0mm 0mm 0mm, clip]{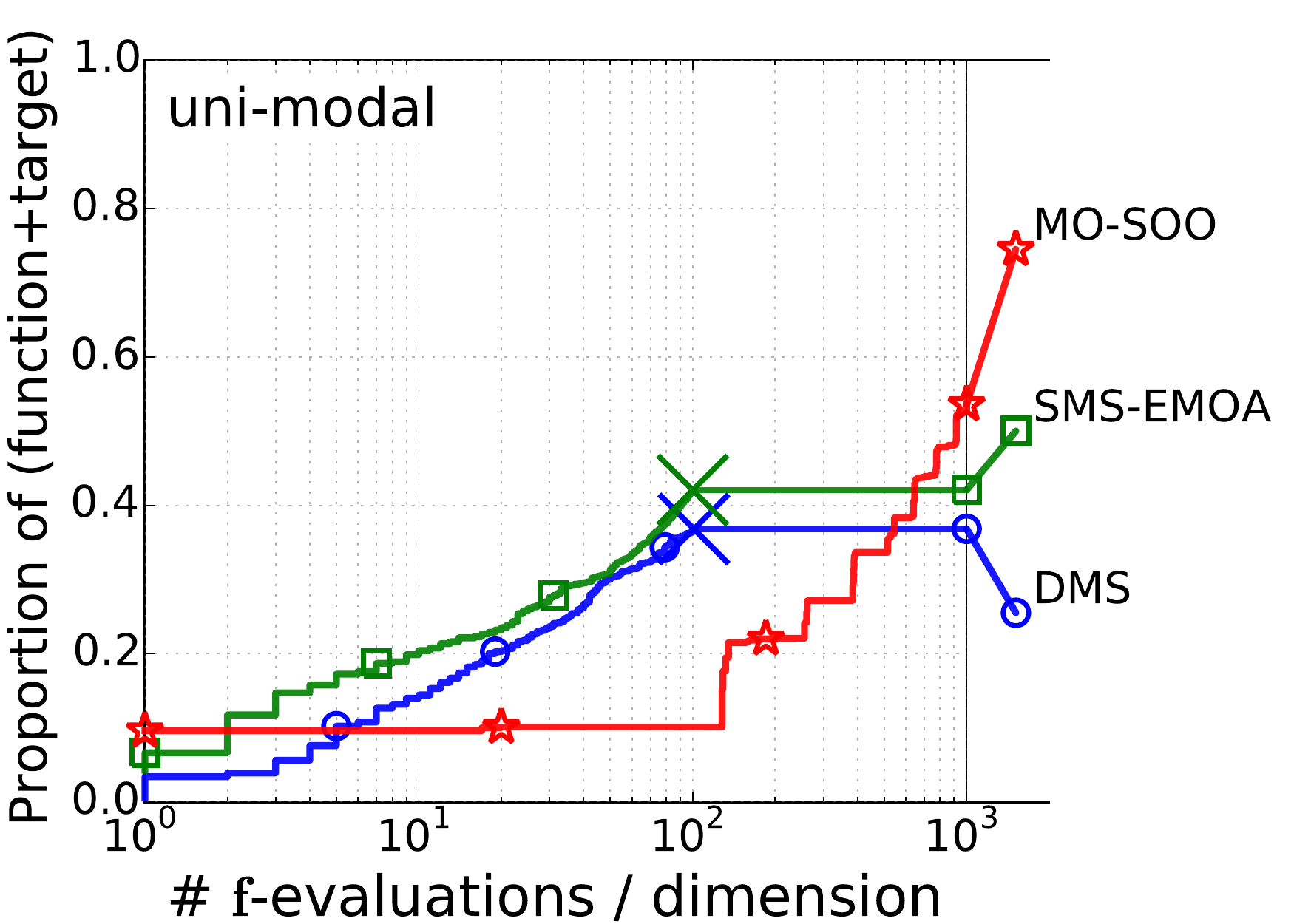}&
				\includegraphics[width=7cm, trim = 0mm 0mm 0mm 0mm, clip]{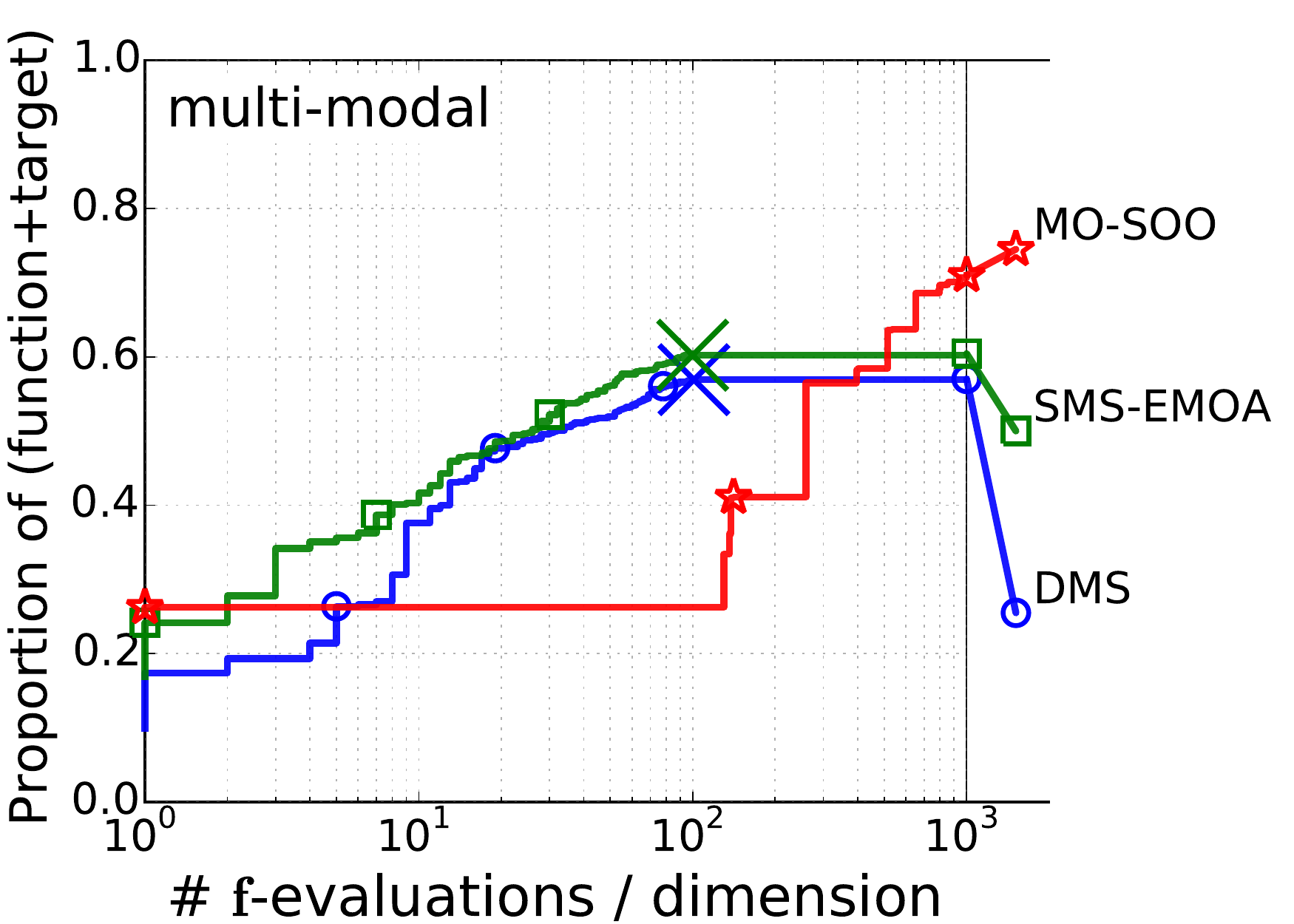}&
				\includegraphics[width=7cm, trim = 0mm 0mm 0mm 0mm, clip]{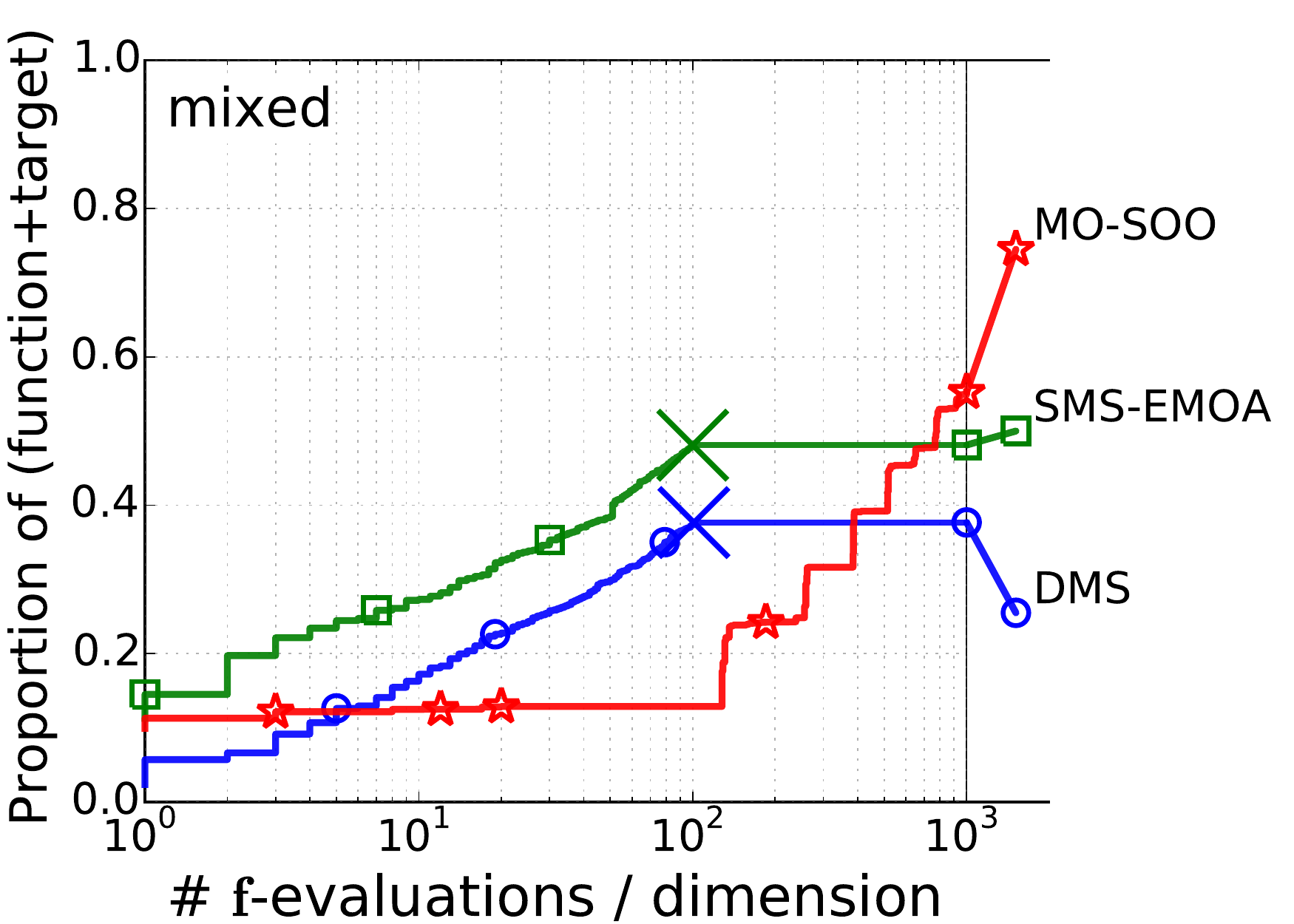}&
				\includegraphics[width=7cm, trim = 0mm 0mm 0mm 0mm, clip]{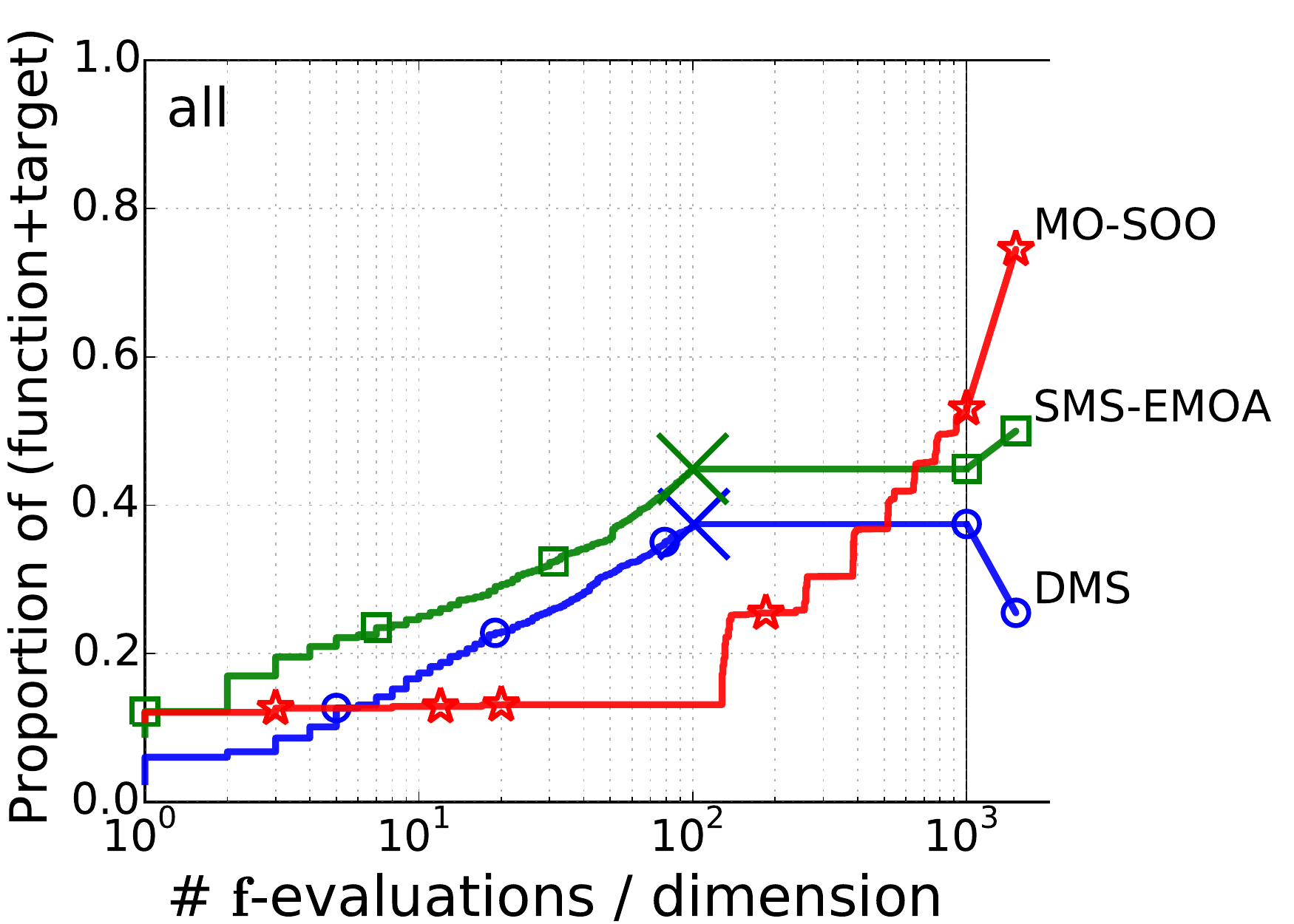}\\
				\toprule
				\multicolumn{4}{c}{\textbf{Additive Epsilon}}\\
				\includegraphics[width=7cm, trim = 0mm 0mm 0mm 0mm, clip]{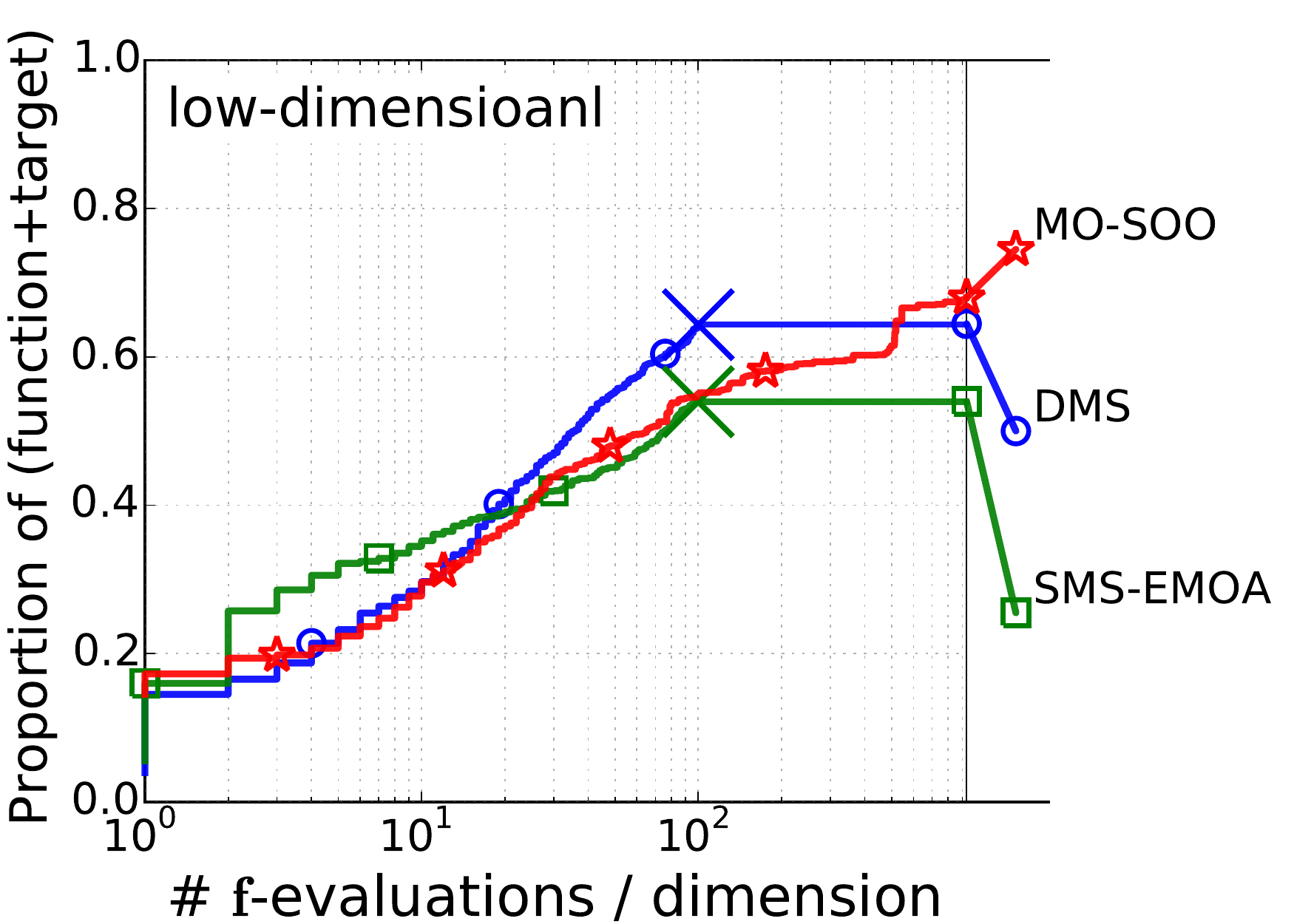}&
				\includegraphics[width=7cm, trim = 0mm 0mm 0mm 0mm, clip]{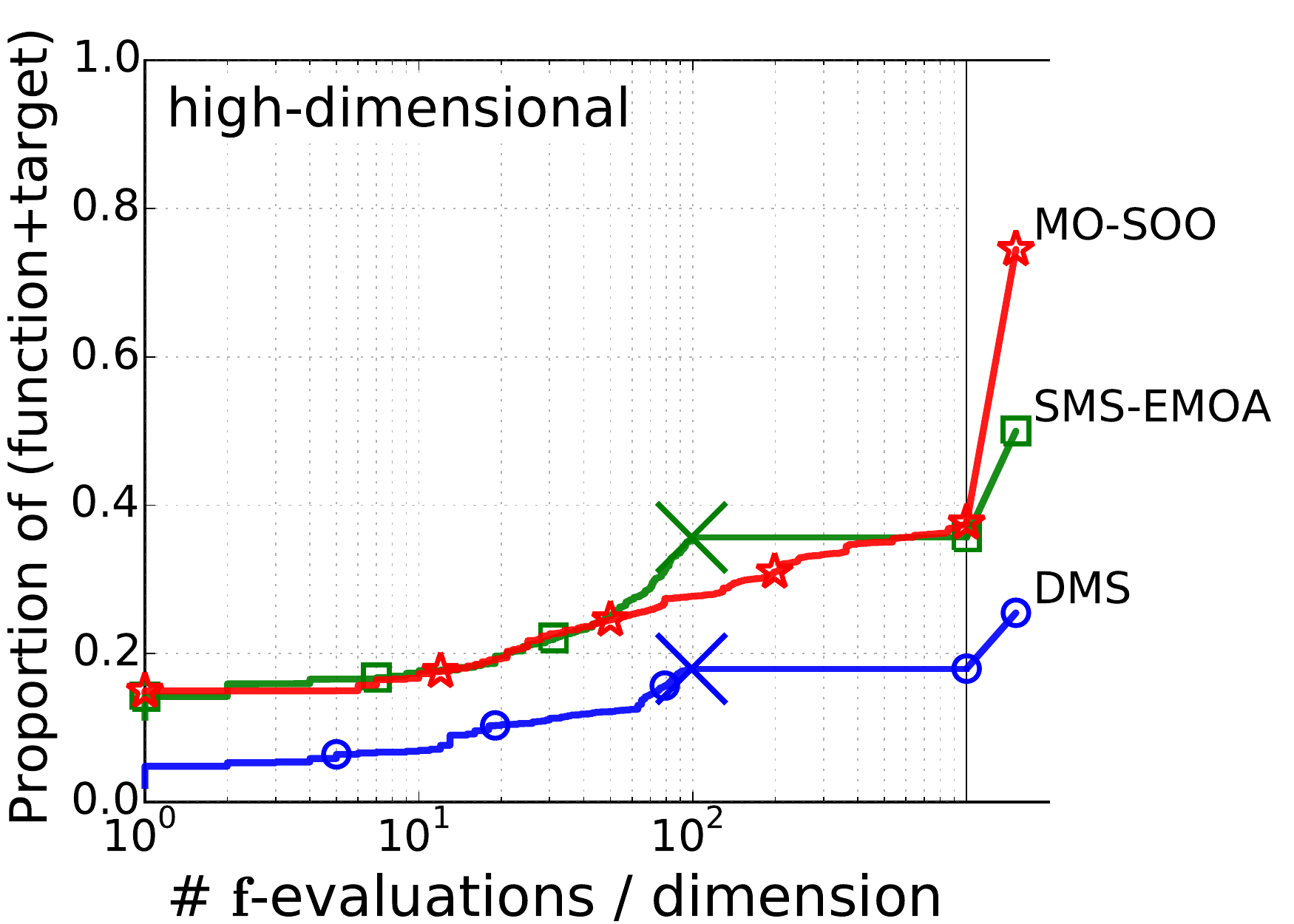}&
				\includegraphics[width=7cm, trim = 0mm 0mm 0mm 0mm, clip]{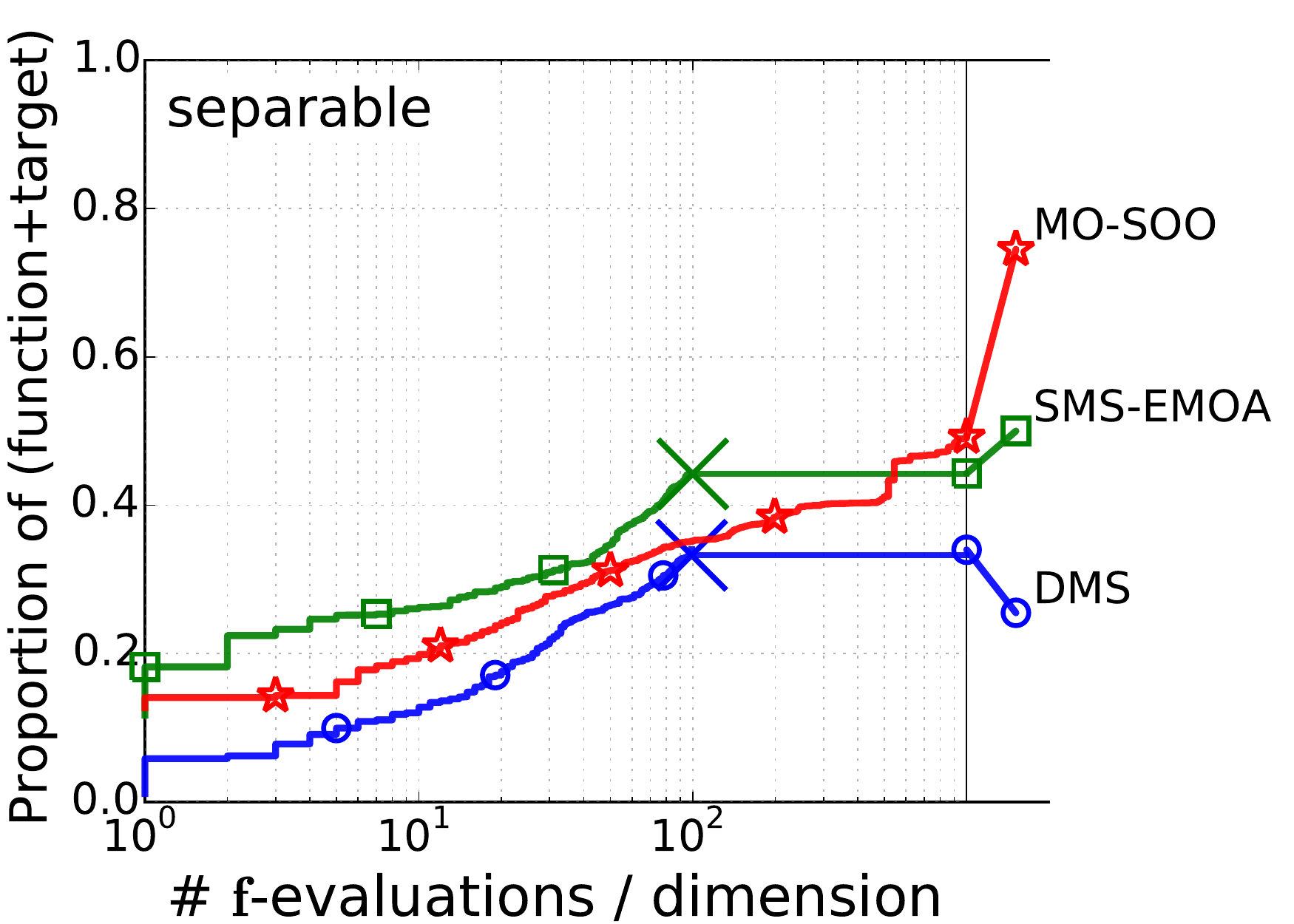}&
				\includegraphics[width=7cm, trim = 0mm 0mm 0mm 0mm, clip]{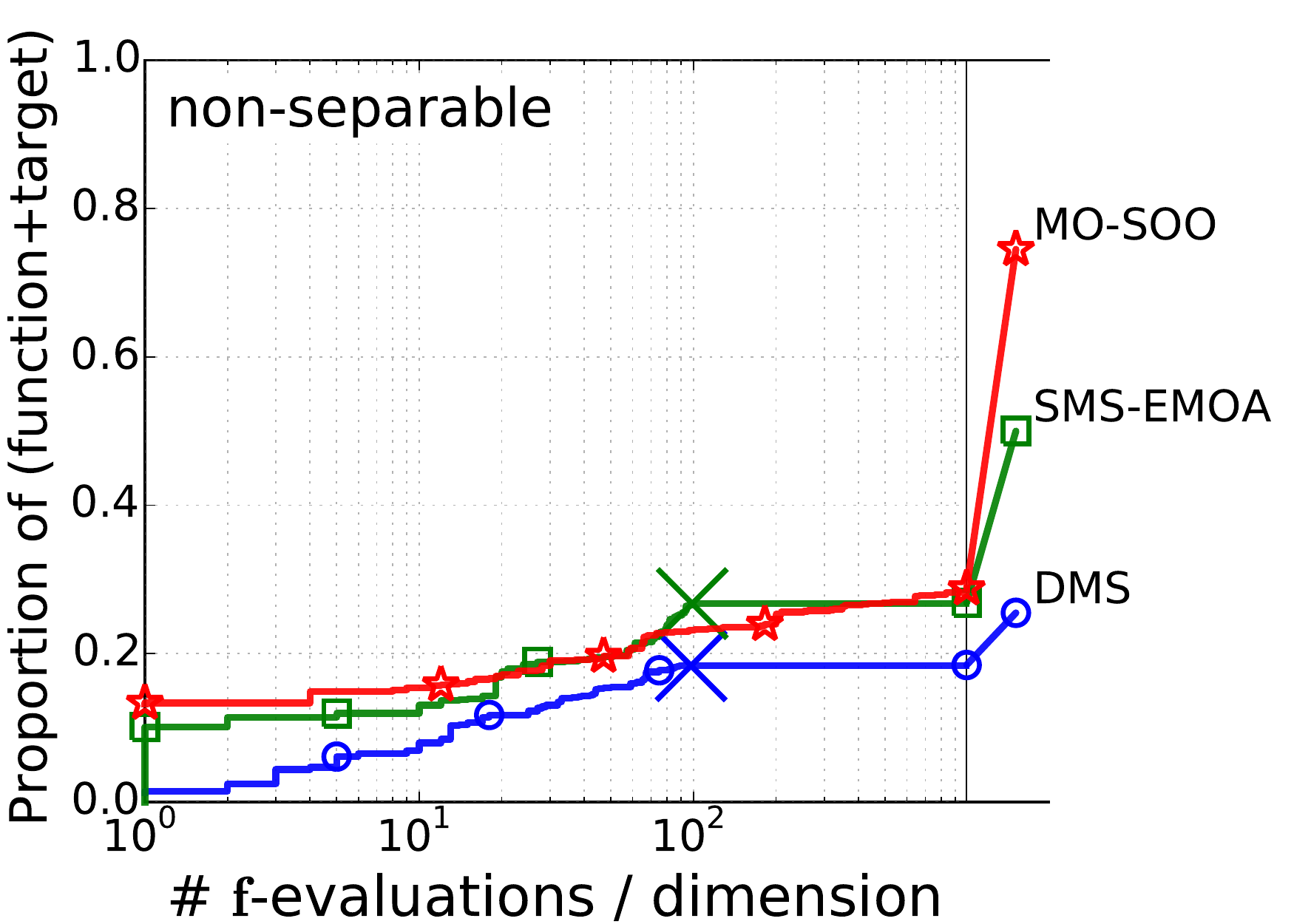}\\
				\includegraphics[width=7cm, trim = 0mm 0mm 0mm 0mm, clip]{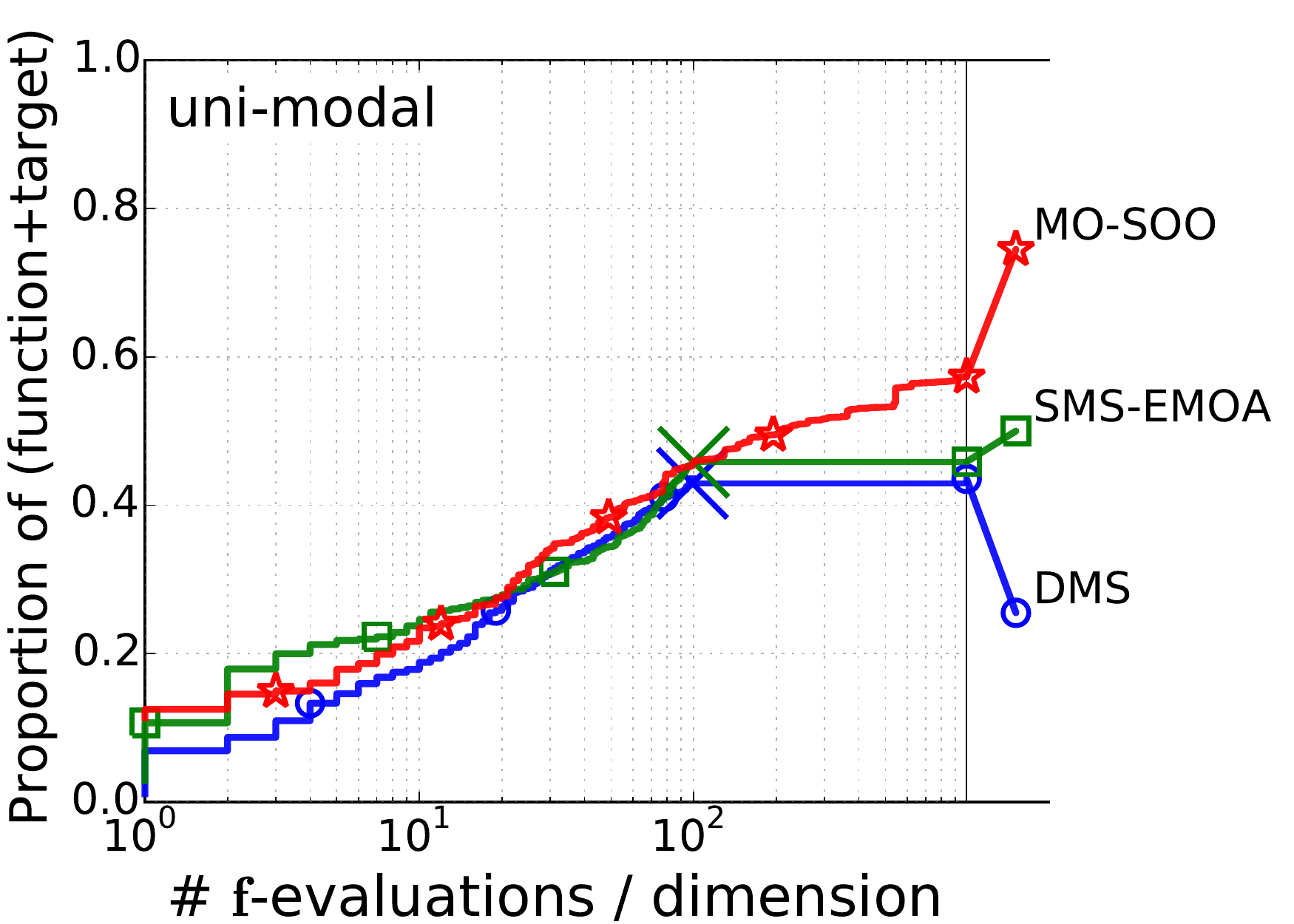}&
				\includegraphics[width=7cm, trim = 0mm 0mm 0mm 0mm, clip]{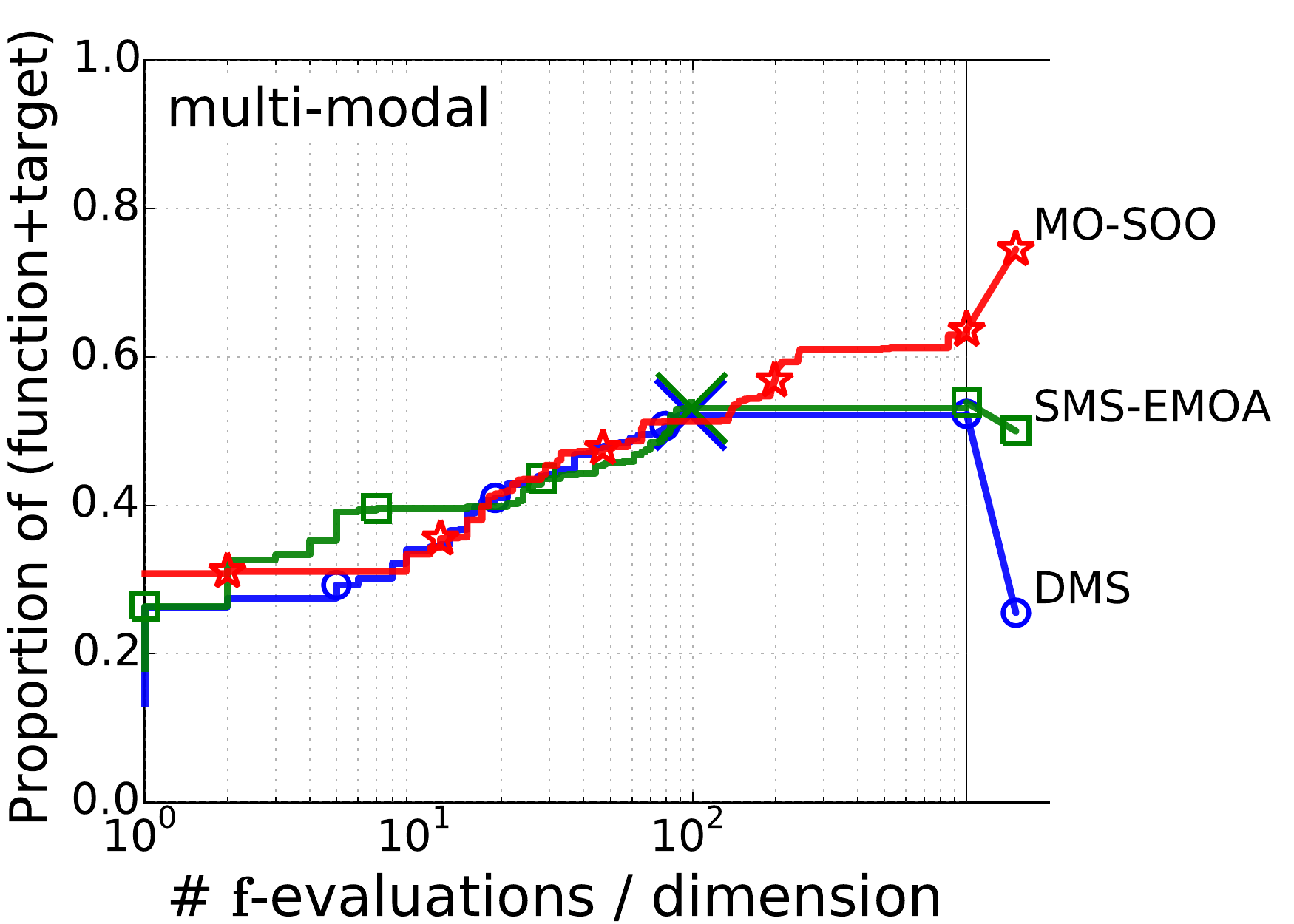}&
					\includegraphics[width=7cm, trim = 0mm 0mm 0mm 0mm, clip]{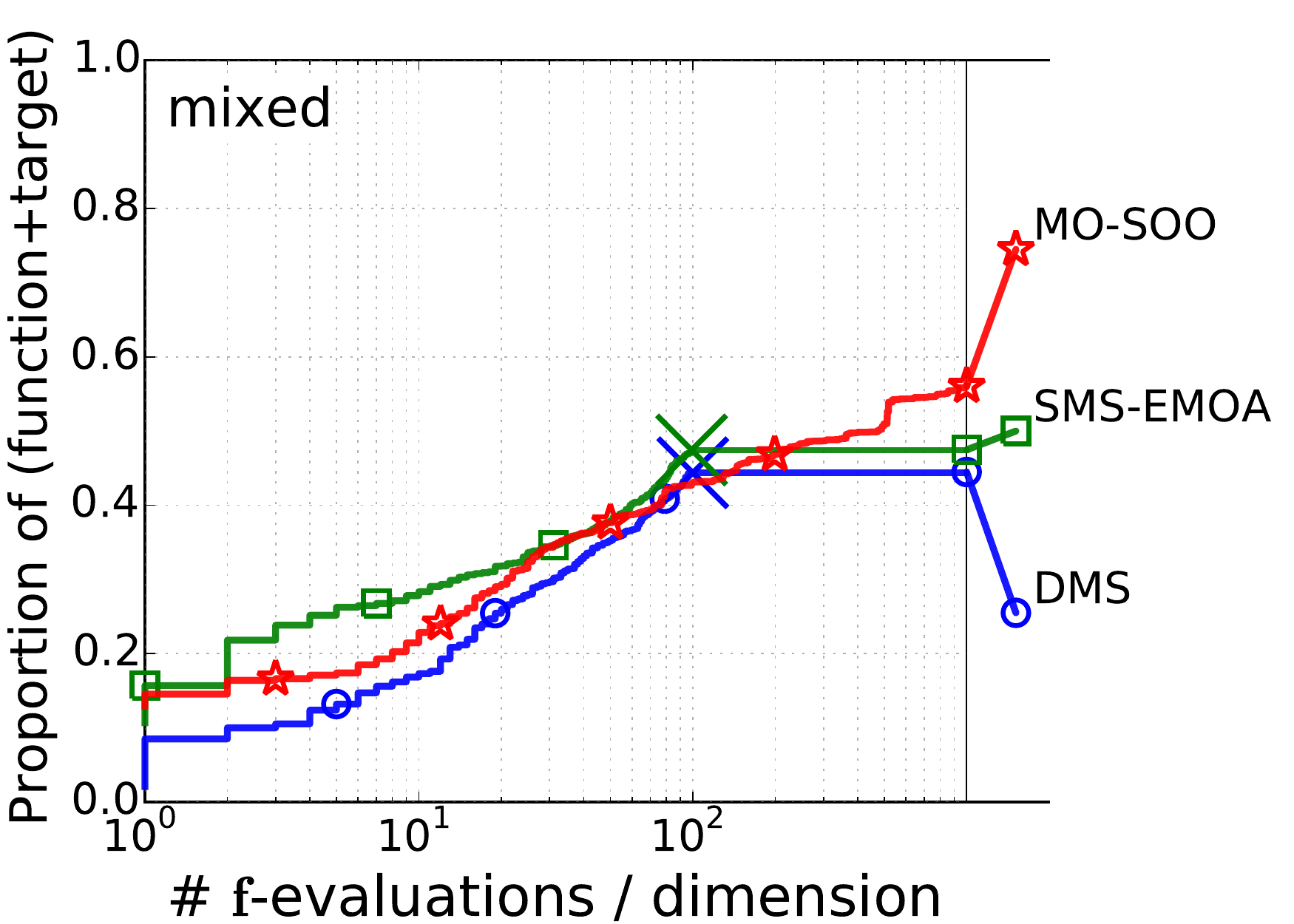}&
					\includegraphics[width=7cm, trim = 0mm 0mm 0mm 0mm, clip]{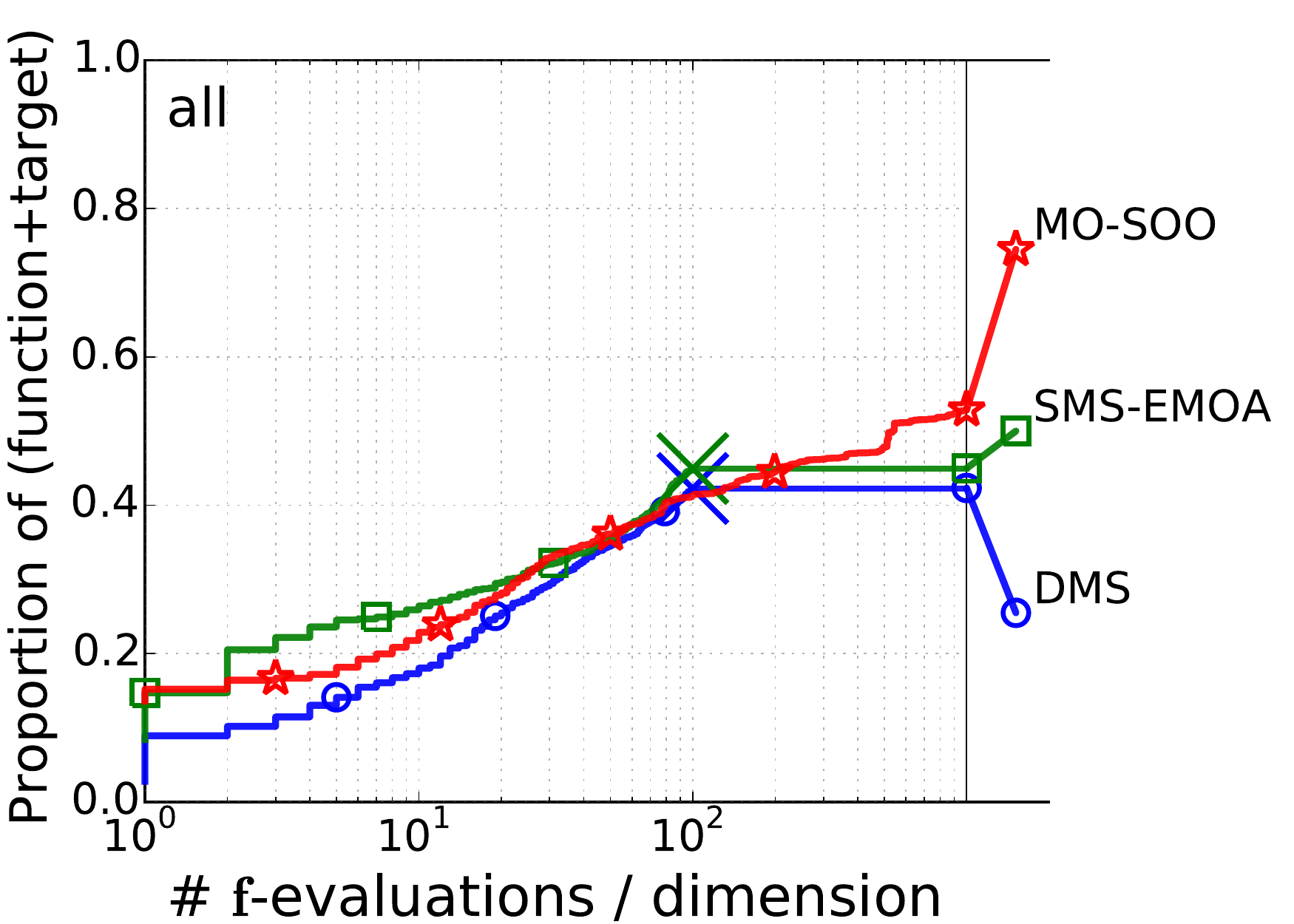}\\
				\toprule
				\multicolumn{4}{c}{\textbf{Generational Distance} (GD)}\\
				\includegraphics[width=7cm, trim = 0mm 0mm 0mm 0mm, clip]{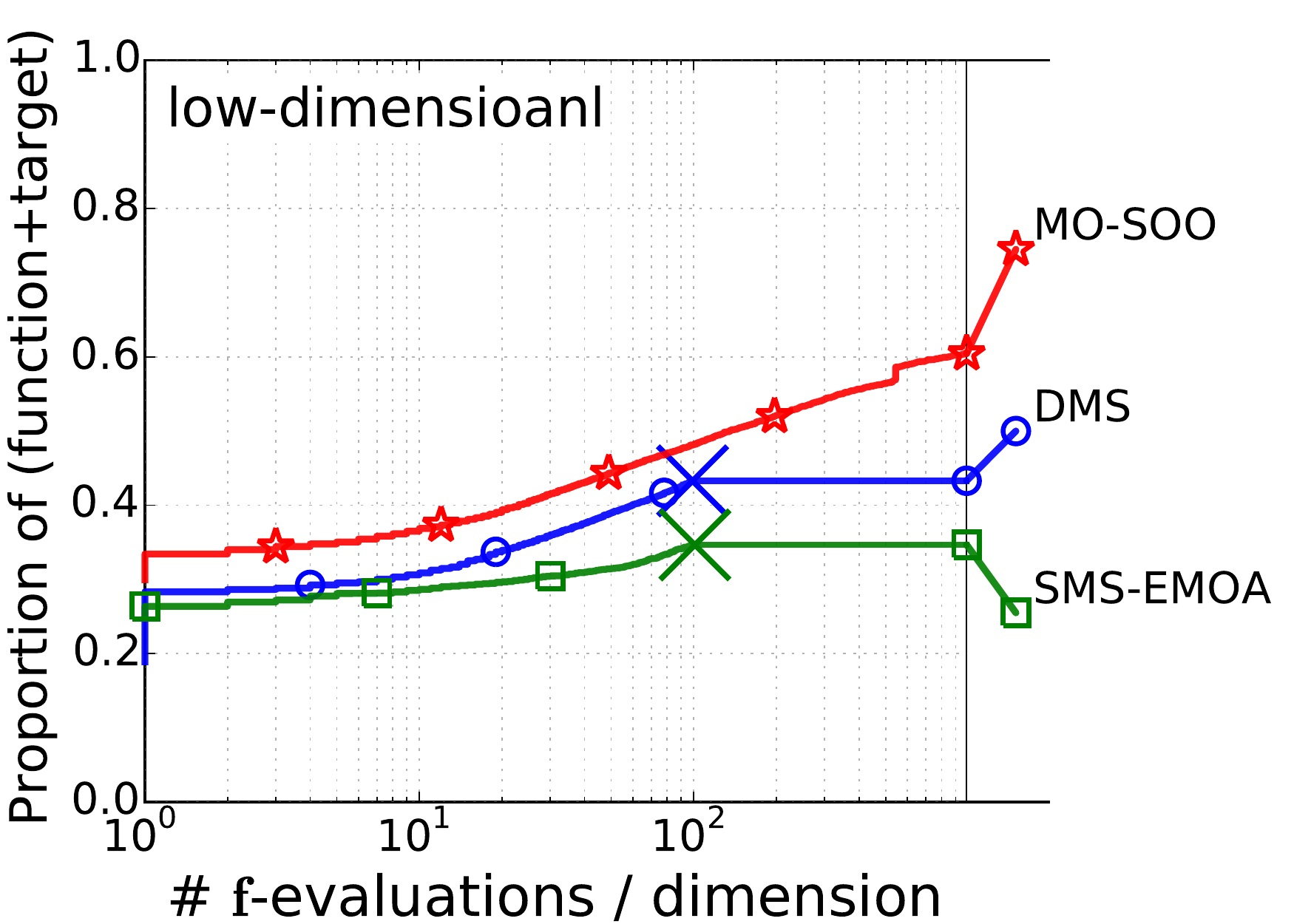}&
				\includegraphics[width=7cm, trim = 0mm 0mm 0mm 0mm, clip]{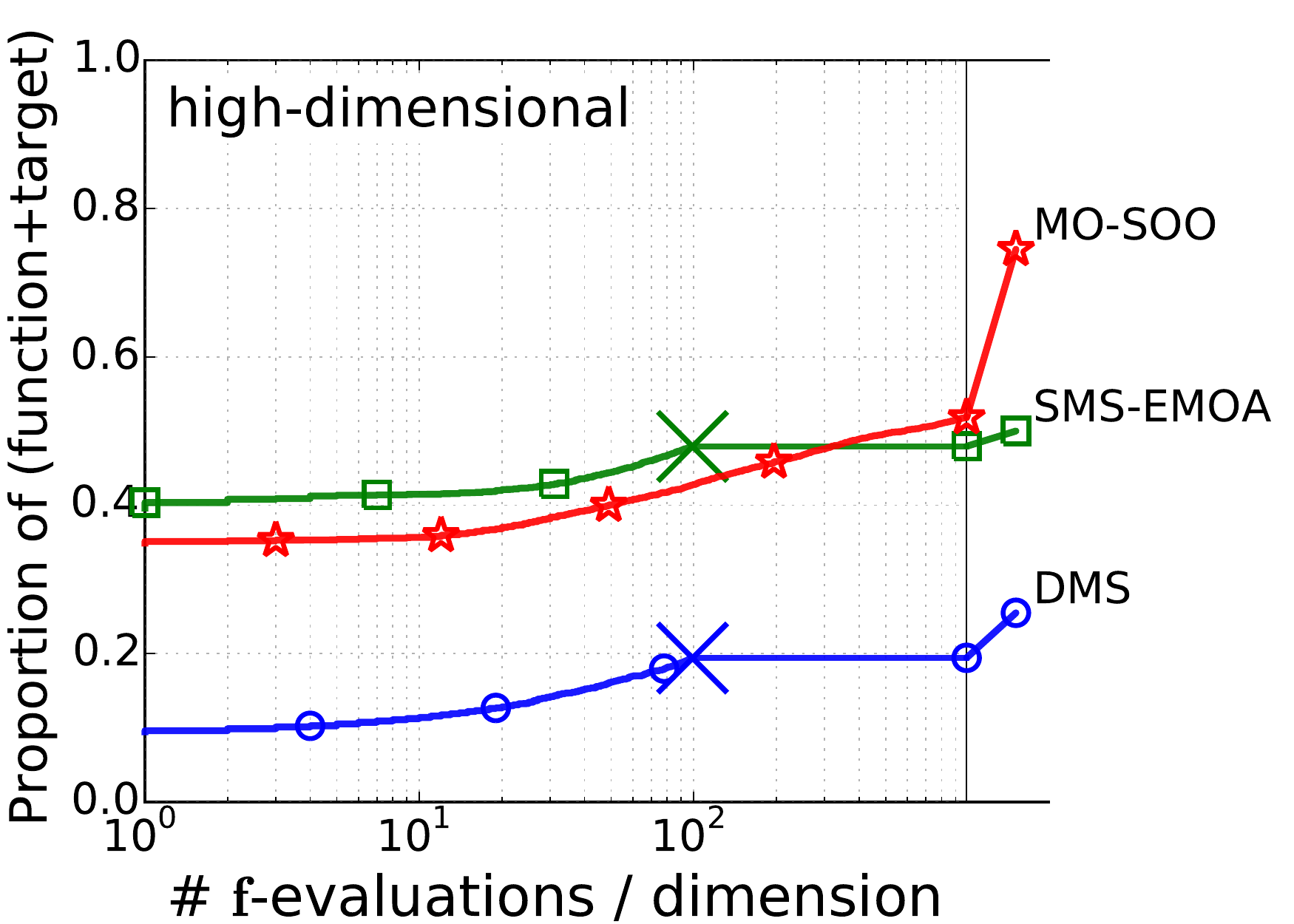}&
				\includegraphics[width=7cm, trim = 0mm 0mm 0mm 0mm, clip]{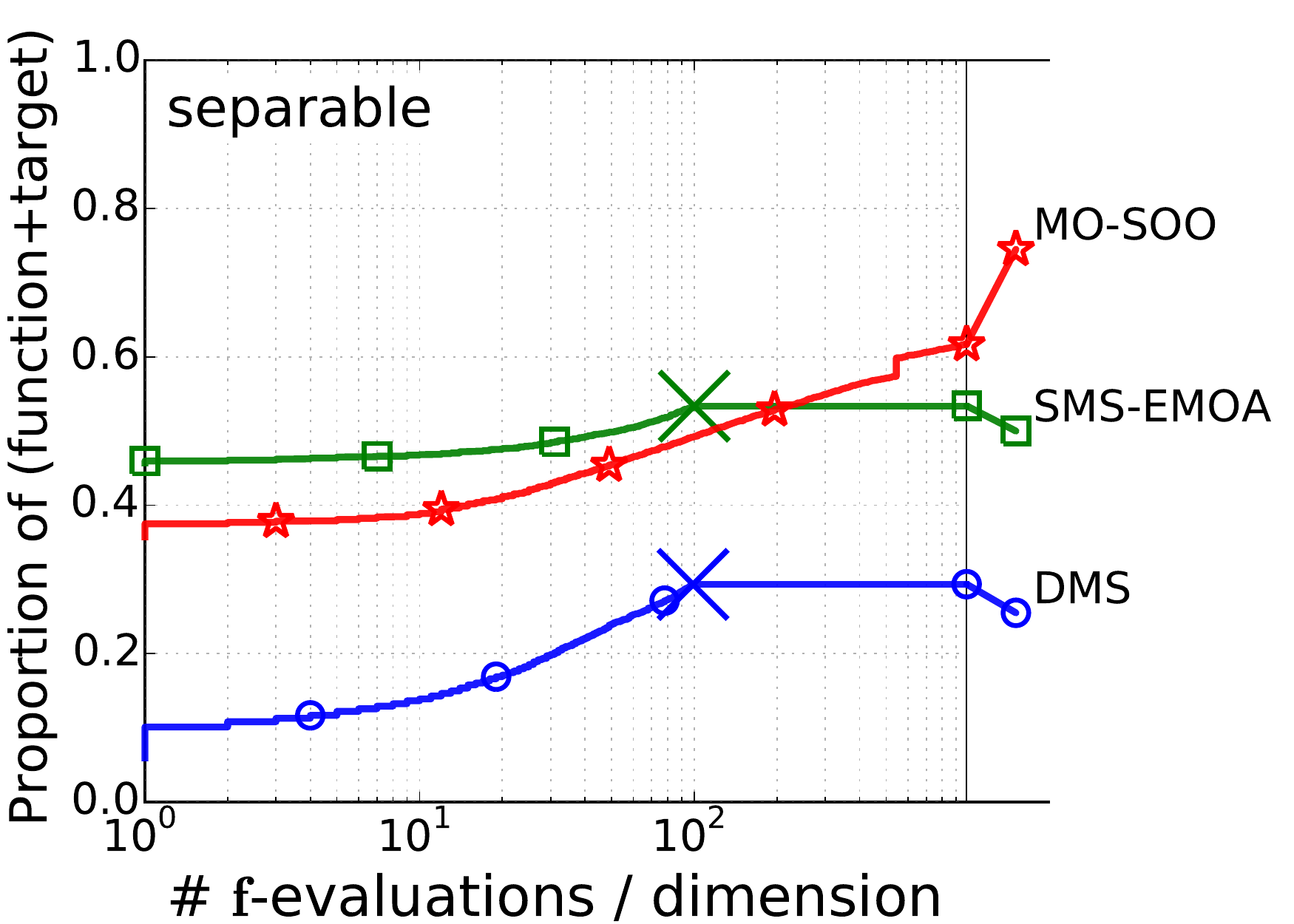}&
				\includegraphics[width=7cm, trim = 0mm 0mm 0mm 0mm, clip]{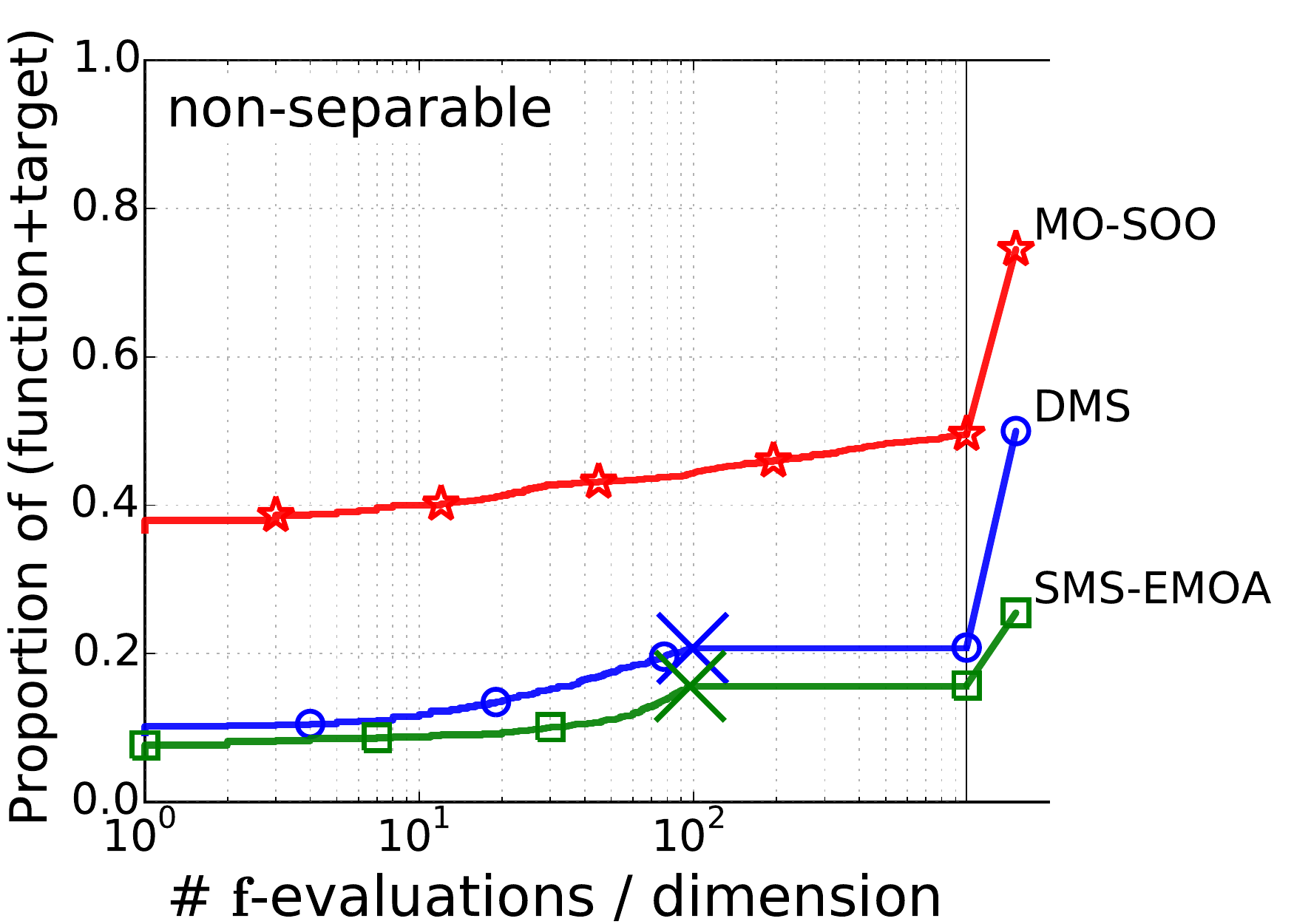}\\
				\includegraphics[width=7cm, trim = 0mm 0mm 0mm 0mm, clip]{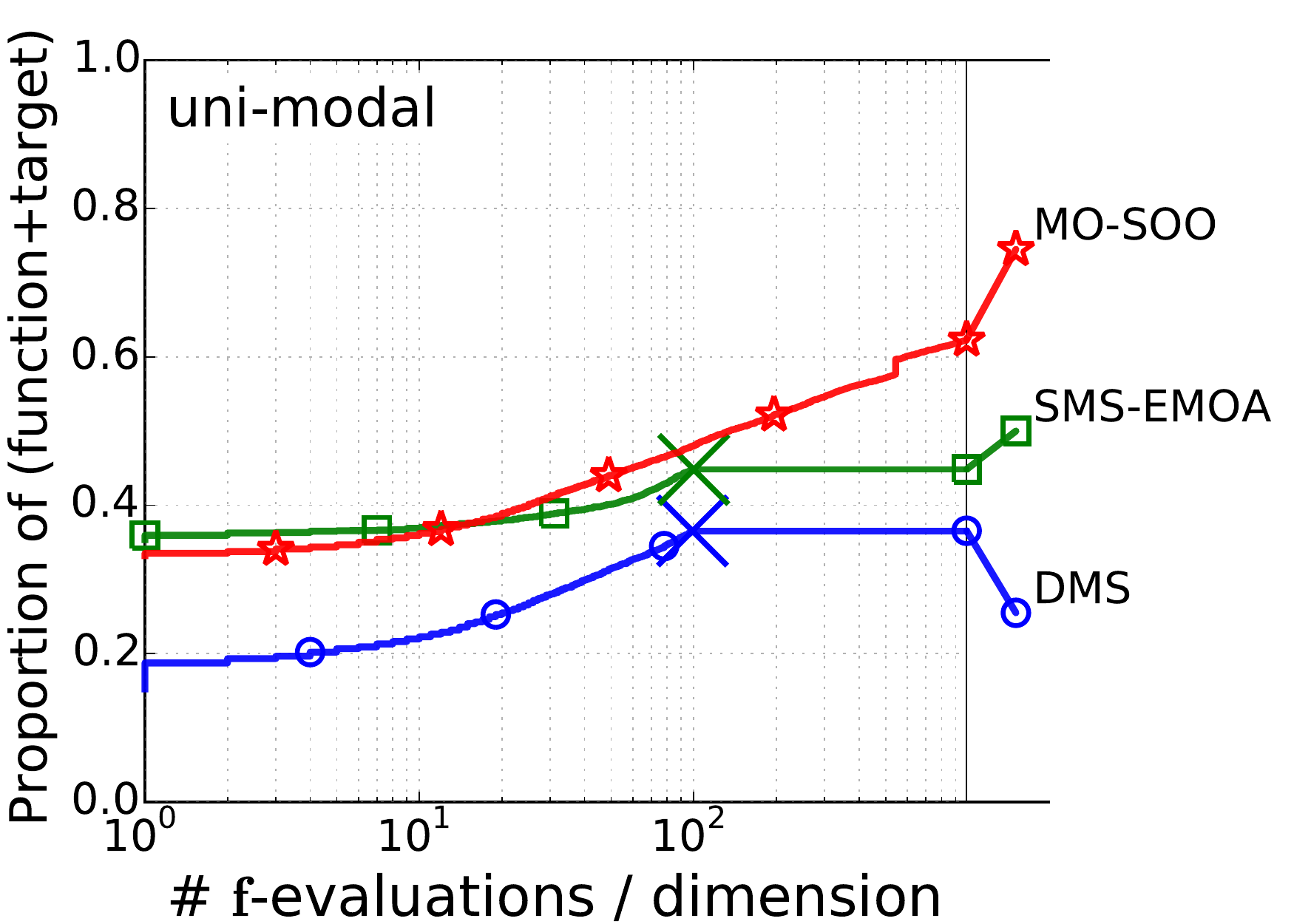}&
				\includegraphics[width=7cm, trim = 0mm 0mm 0mm 0mm, clip]{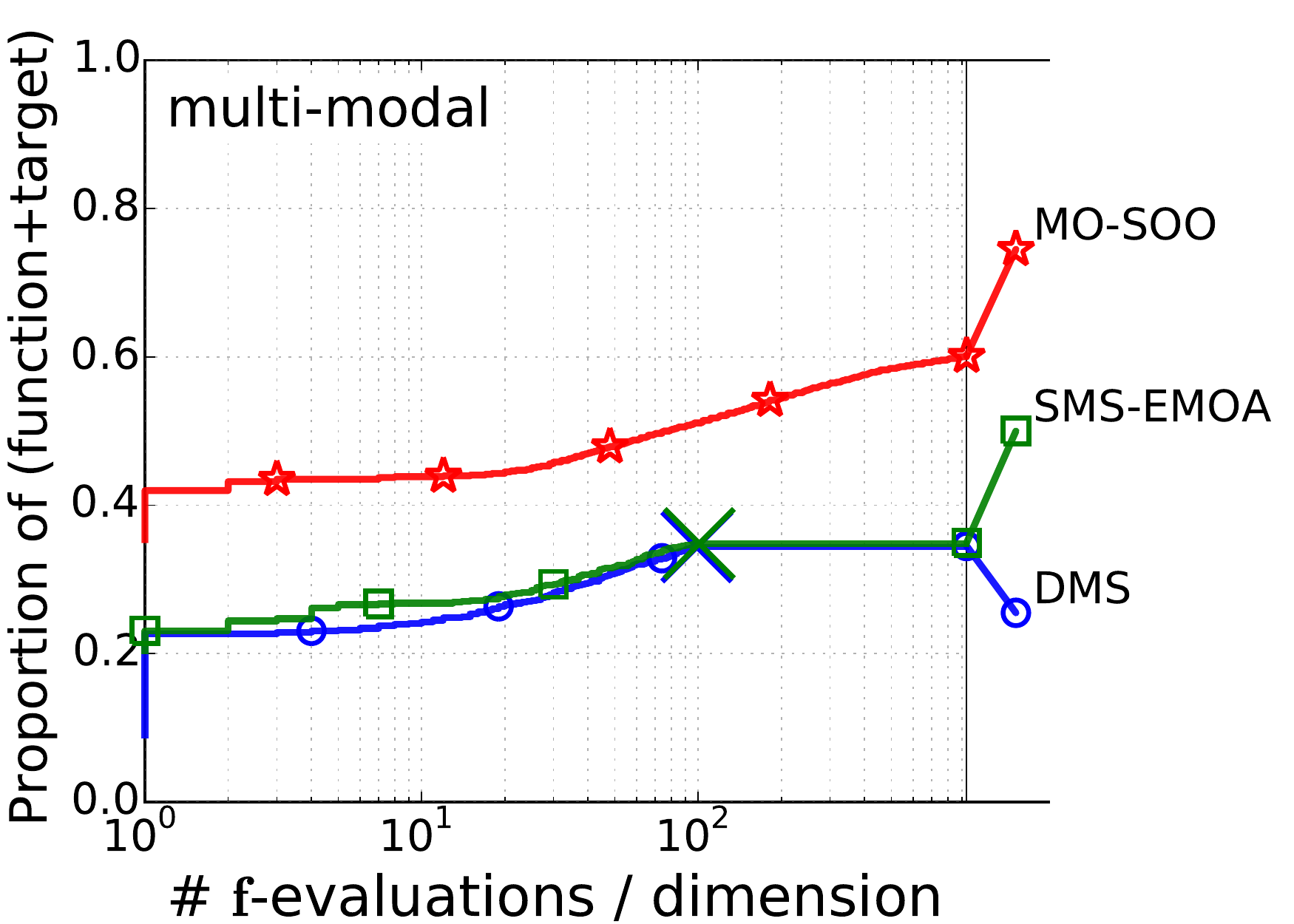}&
				\includegraphics[width=7cm, trim = 0mm 0mm 0mm 0mm, clip]{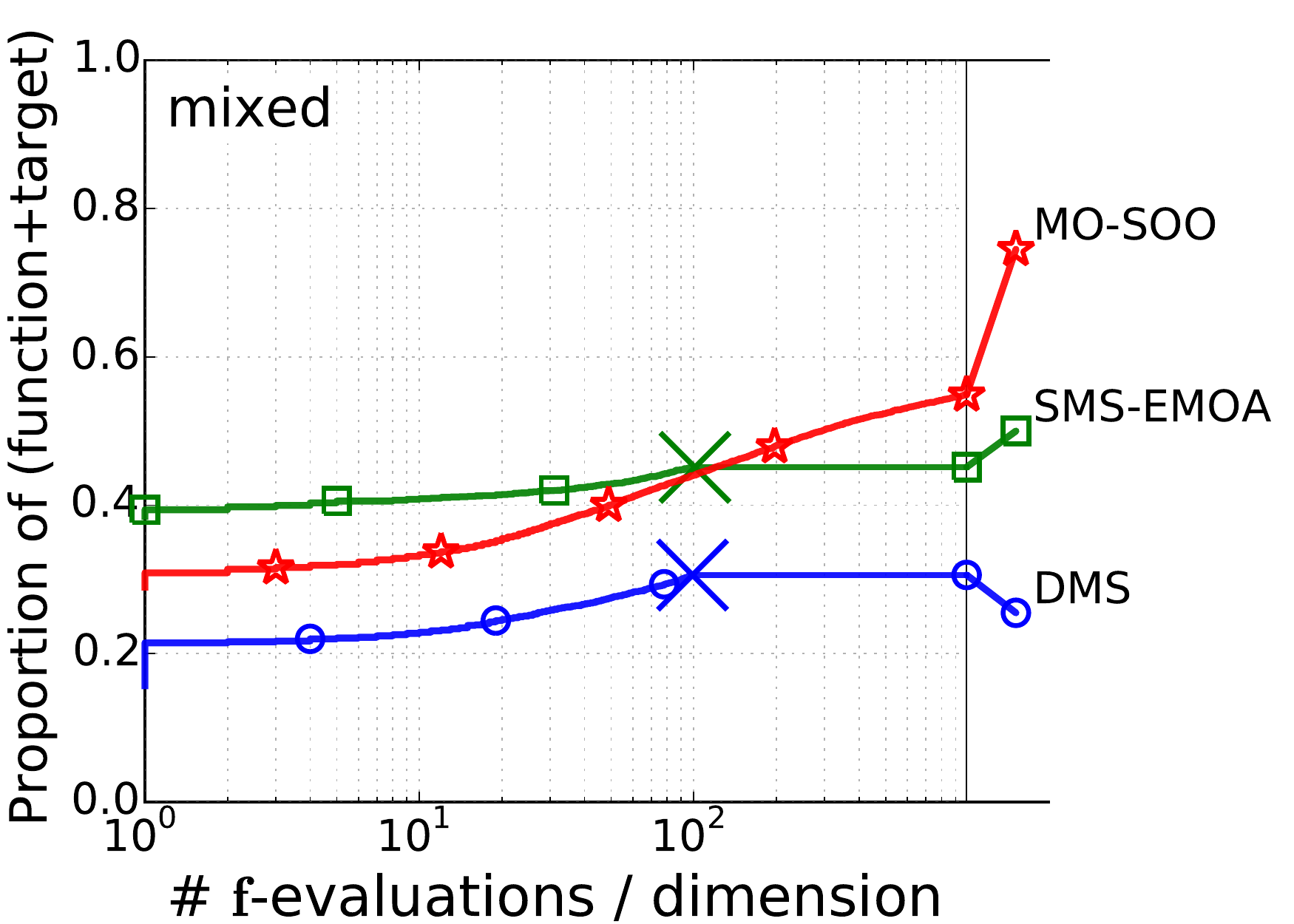}&
				\includegraphics[width=7cm, trim = 0mm 0mm 0mm 0mm, clip]{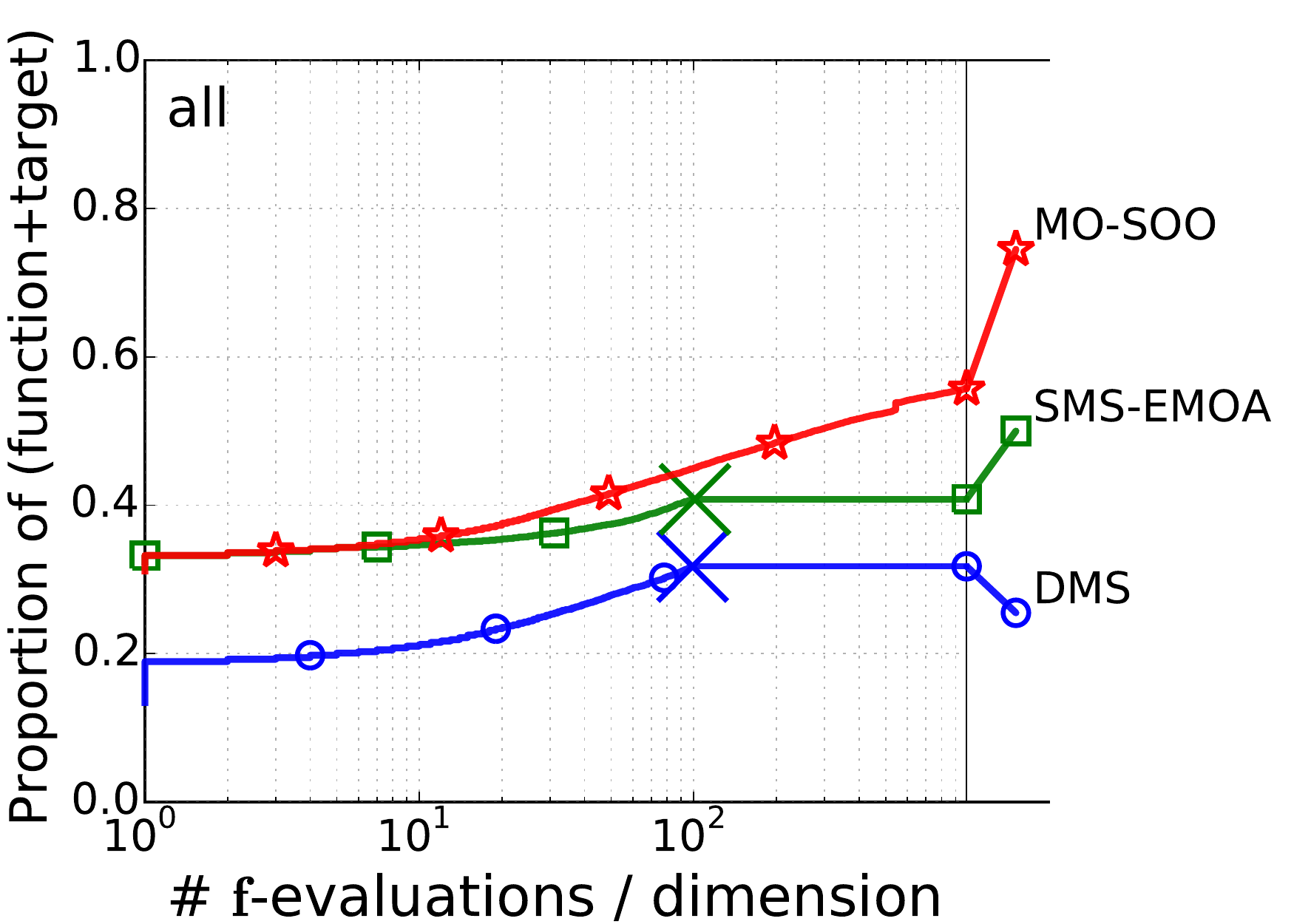}\\
				\toprule
				\multicolumn{4}{c}{\textbf{Inverted Generational Distance} (IGD)}\\
				\includegraphics[width=7cm, trim = 0mm 0mm 0mm 0mm, clip]{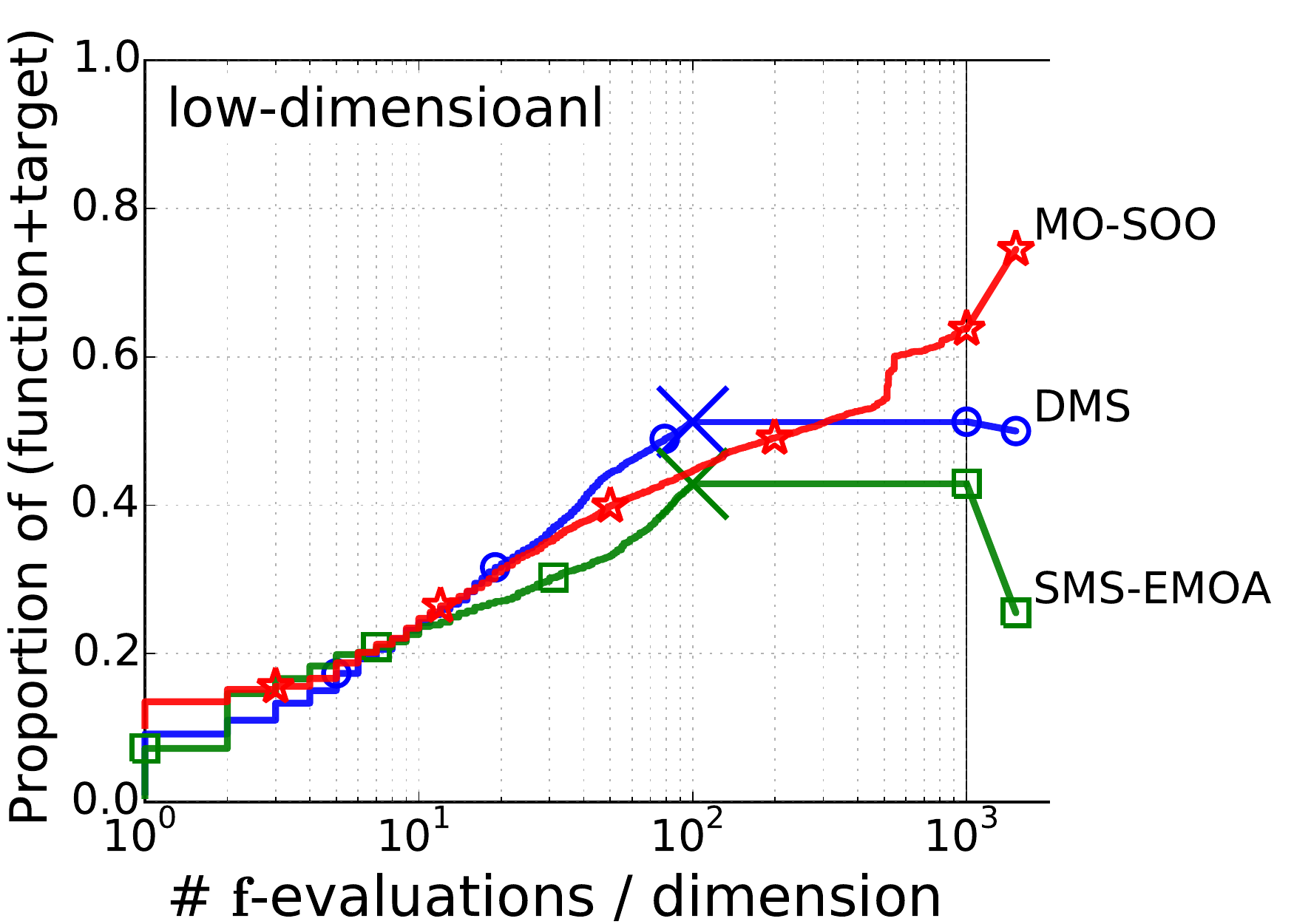}&
				\includegraphics[width=7cm, trim = 0mm 0mm 0mm 0mm, clip]{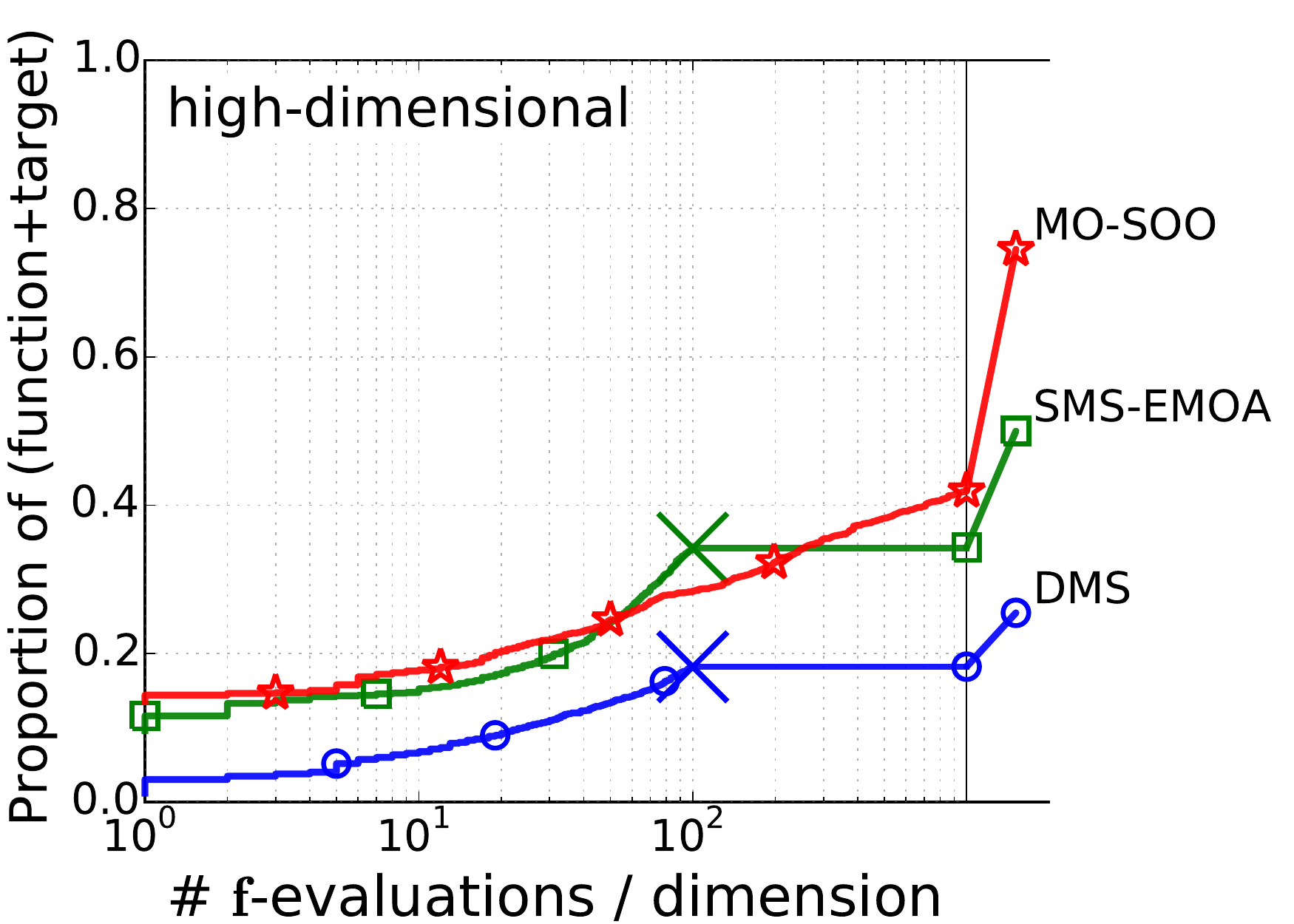}&
				\includegraphics[width=7cm, trim = 0mm 0mm 0mm 0mm, clip]{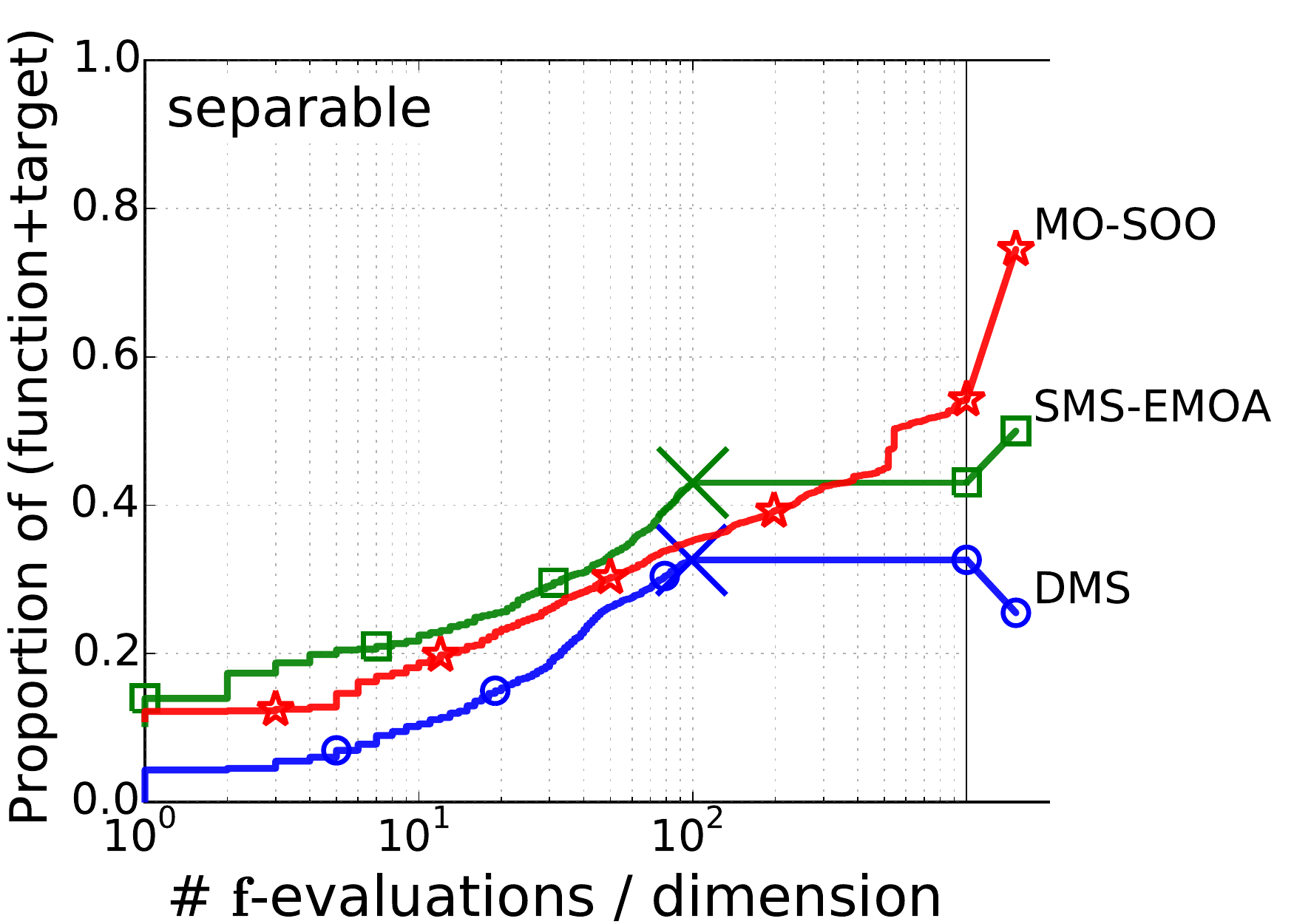}&
				\includegraphics[width=7cm, trim = 0mm 0mm 0mm 0mm, clip]{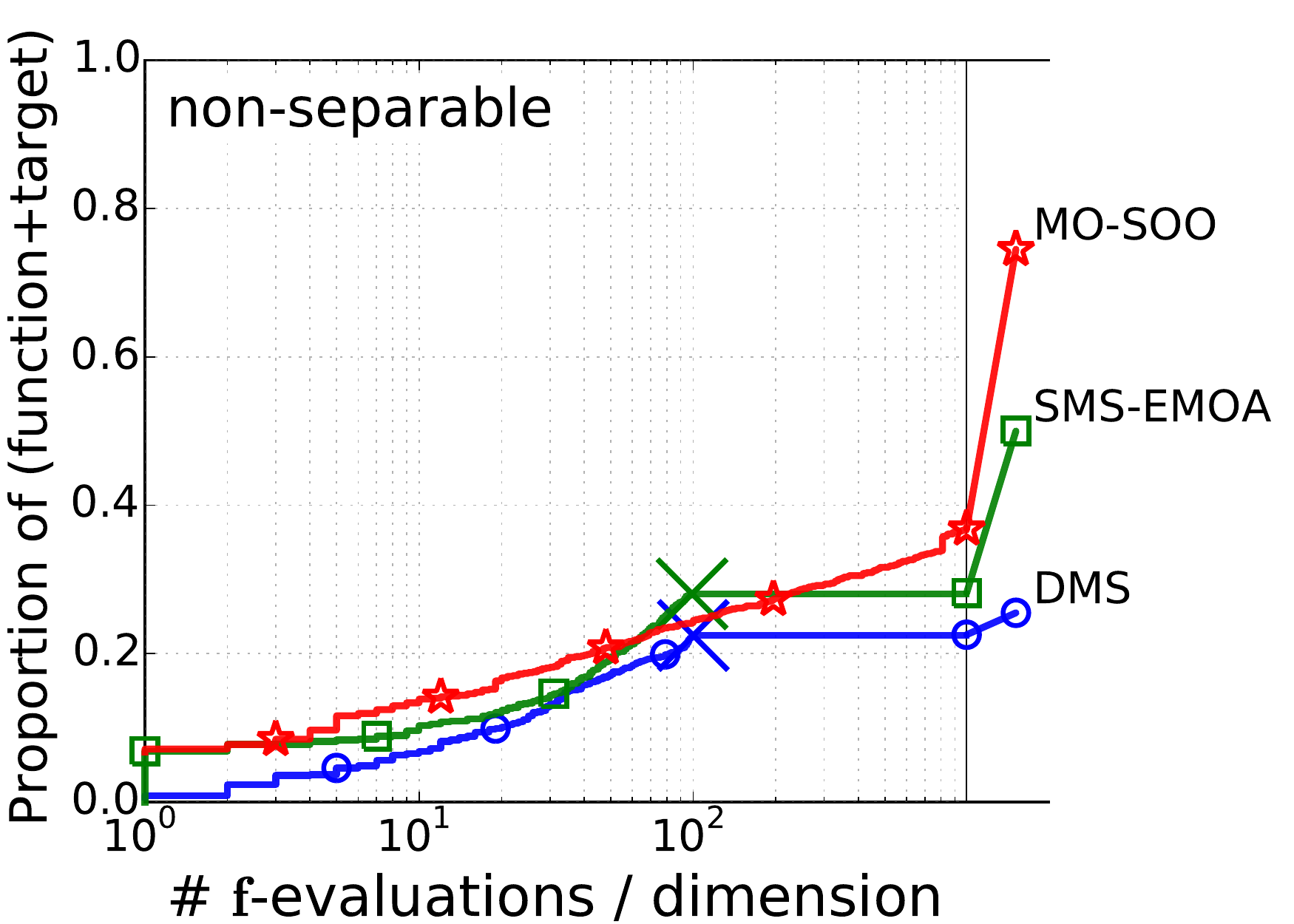}\\
				\includegraphics[width=7cm, trim = 0mm 0mm 0mm 0mm, clip]{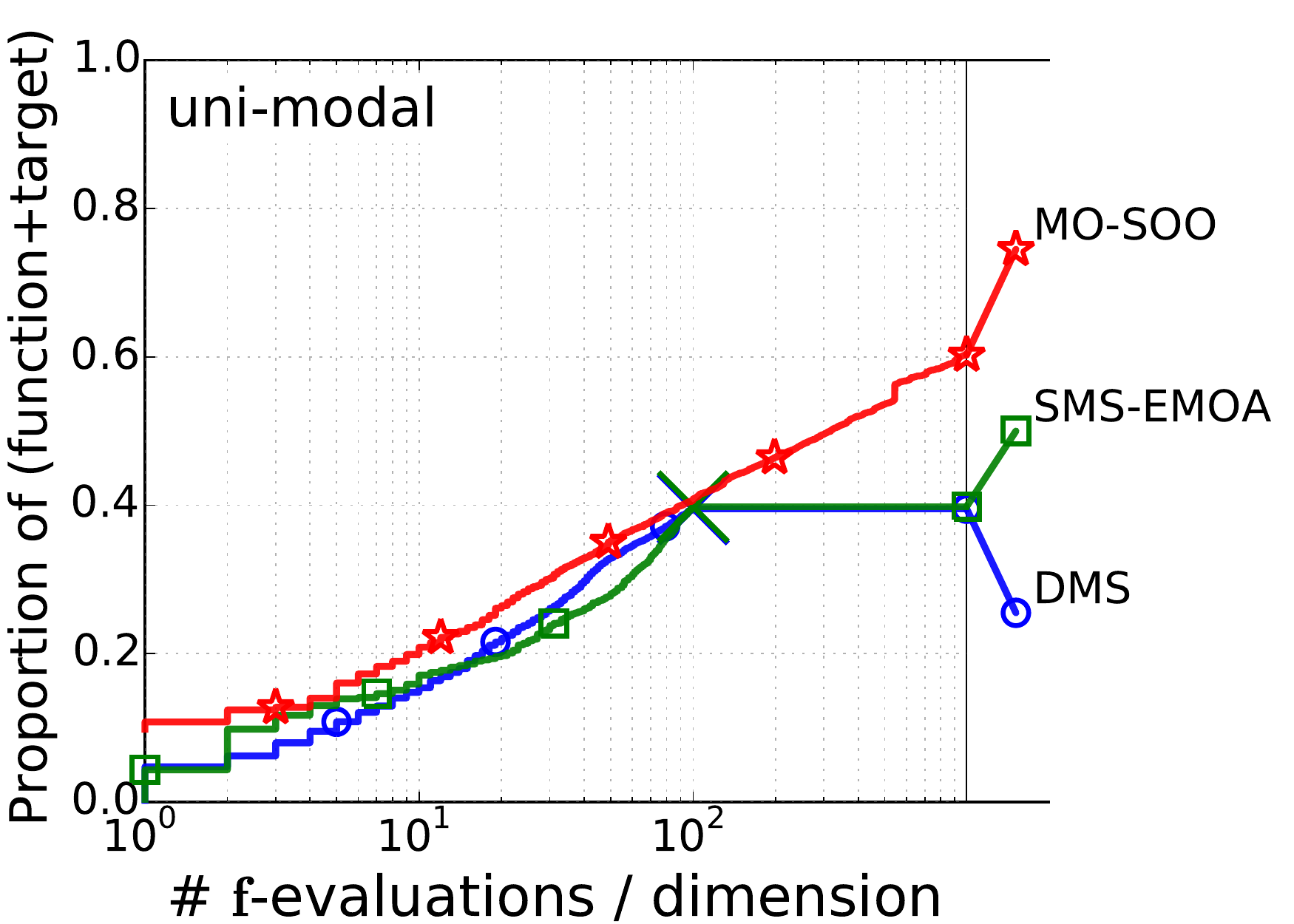}&
				\includegraphics[width=7cm, trim = 0mm 0mm 0mm 0mm, clip]{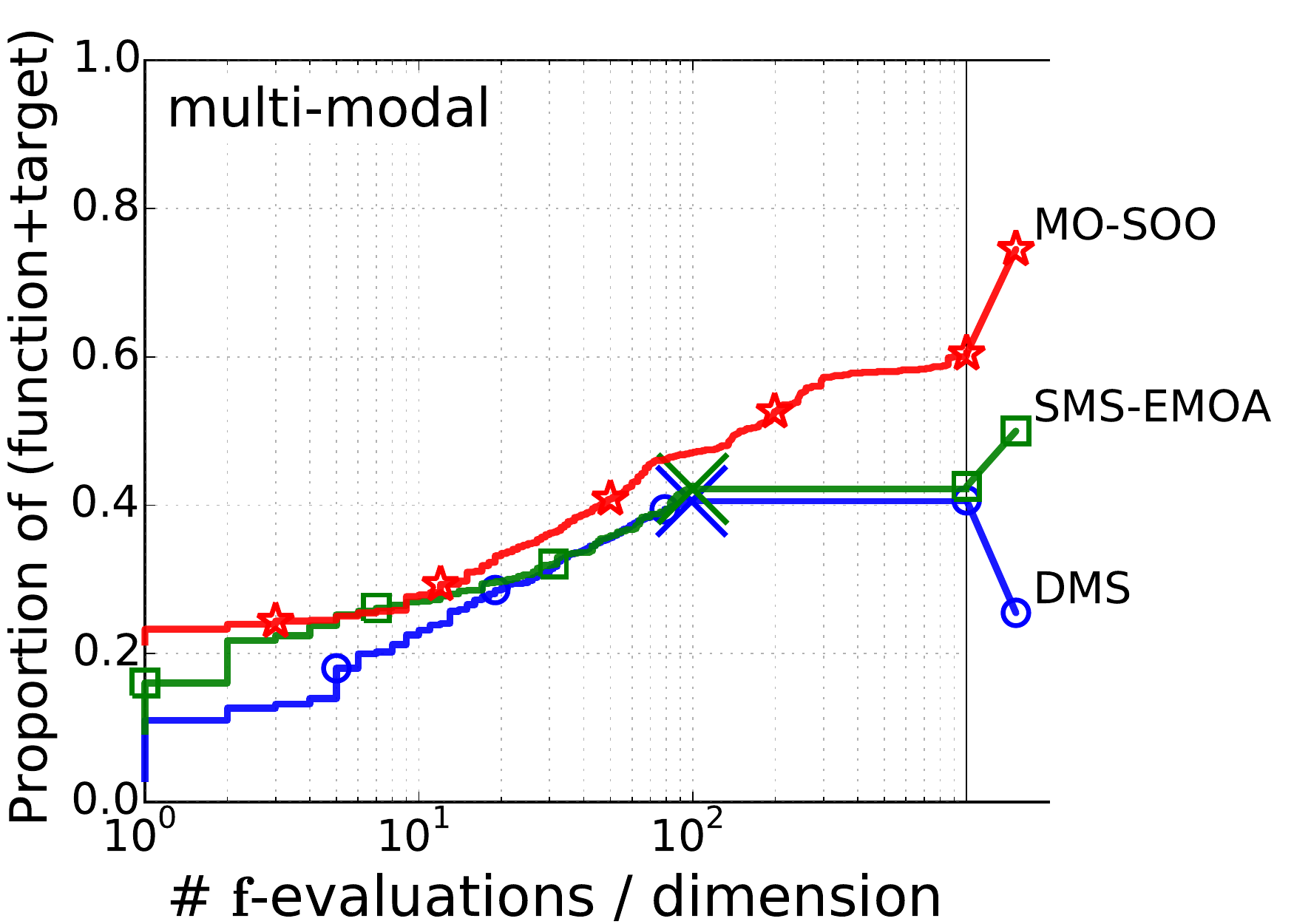}&
				\includegraphics[width=7cm, trim = 0mm 0mm 0mm 0mm, clip]{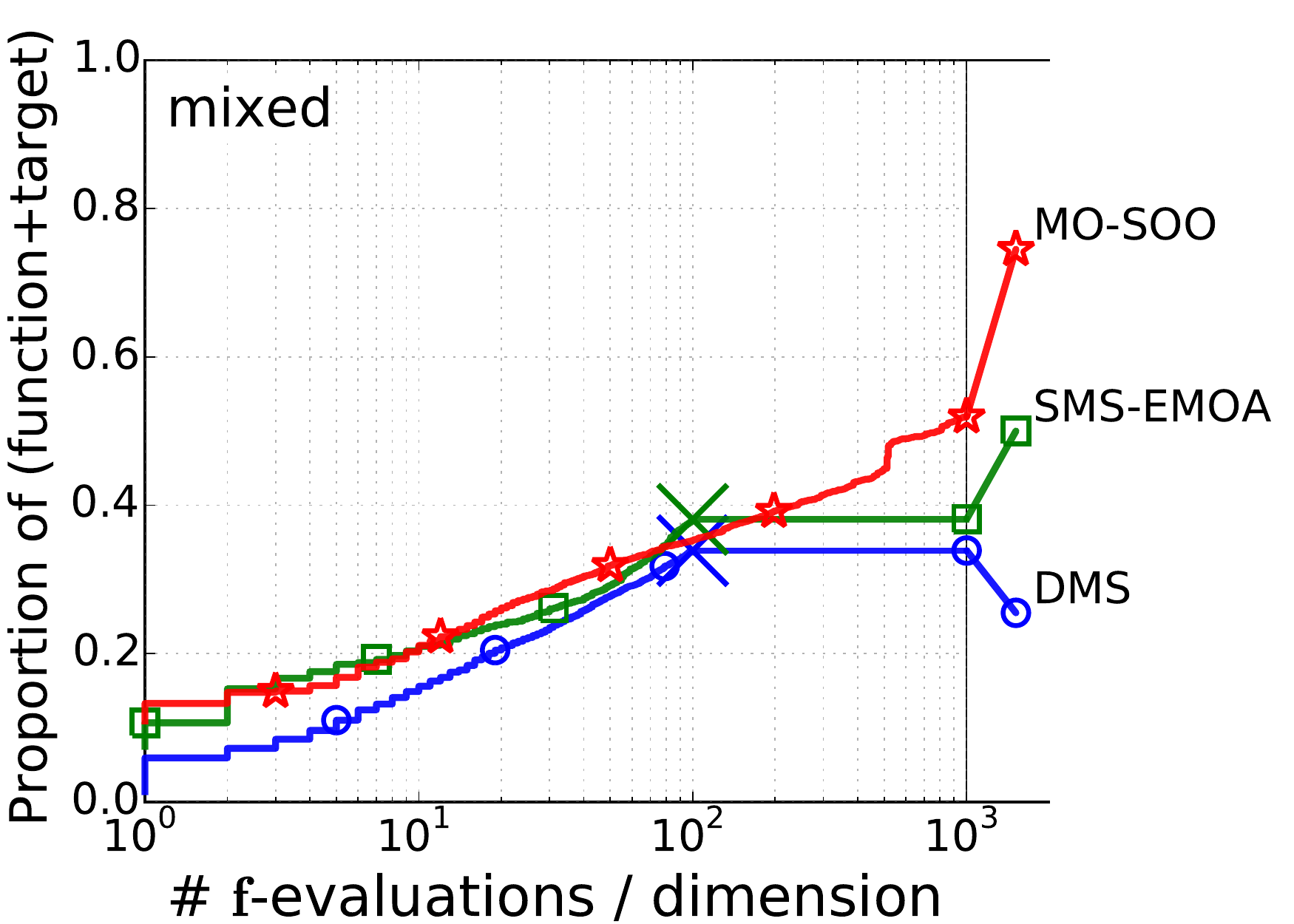}&
				\includegraphics[width=7cm, trim = 0mm 0mm 0mm 0mm, clip]{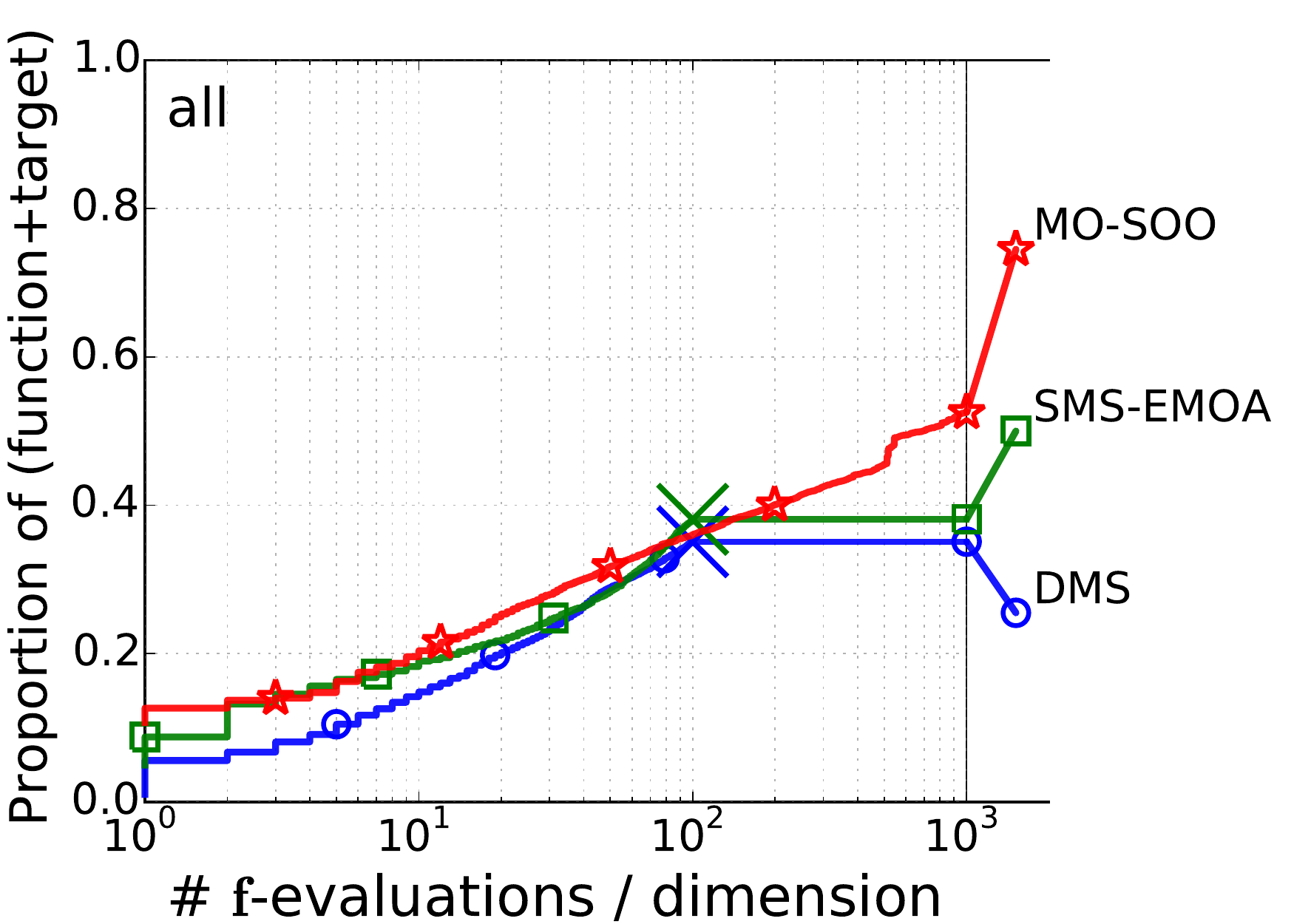}\\
				\bottomrule
			\end{tabular}}
		\end{center}
		\caption{Data profiles aggregated over problem categories for each of the quality indicators computed. The symbol $\times$ indicates the maximum number of function
			evaluations.}
		\label{fig:hv_agg_per_dim}		
	\end{figure*}

Figures~\ref{fig:all_ind} and~\ref{fig:hv_agg_per_dim} show the data profiles of the compared algorithms as a function of the number of functione evaluations used.

			\section{Empirical Runtime Evaluation}
			\label{sec:timing}
			
			\begin{figure}[tb]
				\centering
				\includegraphics[width=1\textwidth]{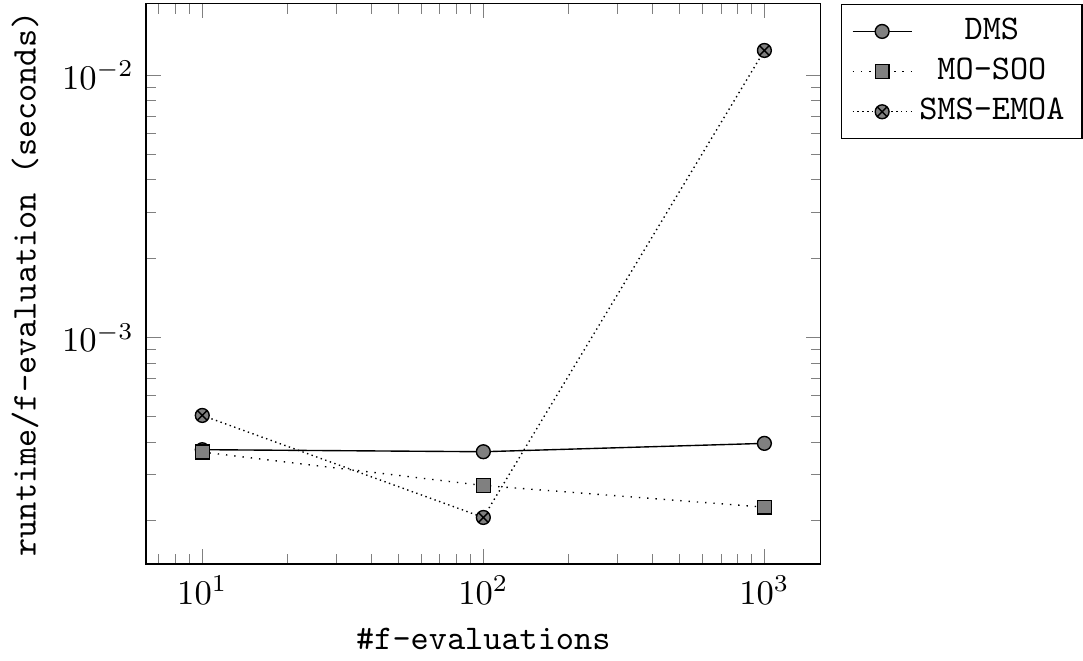}
				\caption{A $log$-$\log$ plot visualizing the runtime per one function evaluation (in seconds) of the compared algorithms. All the algorithms were run on a selected set of problems over a set of evaluation budgets, namely BK1, DPAM1, L3ZDT1, DTLZ3, and FES3; with an evaluation budget $\in \{10, 100, 1000\}$
					per problem on a PC with: 64-bit Windows 7, { Intel Xeon E5 CPU @ 3.20GHz,  16GB of memory.}}
				\label{fig:timing}
			\end{figure}
			
			In order to evaluate the complexity of the algorithms (measured in runtime), we have run the algorithms on a representative set of the problems.  The empirical complexity of an algorithm is then computed as the running time (in seconds) of the algorithm summed over all the problems given an evaluation budget~($\text{\#FE}$). The results are shown in Figure~\ref{fig:timing}. 
			

\bibliographystyle{plainnat}
\bibliography{bmobench}

\begin{thebibliography}{23}
\providecommand{\natexlab}[1]{#1}
\providecommand{\url}[1]{\texttt{#1}}
\expandafter\ifx\csname urlstyle\endcsname\relax
  \providecommand{\doi}[1]{doi: #1}\else
  \providecommand{\doi}{doi: \begingroup \urlstyle{rm}\Url}\fi

\bibitem[Al-Dujaili and Suresh(2017)]{ash-mosoo-15}
Abdullah Al-Dujaili and S.~Suresh.
\newblock Multi-objective simultaneous optimistic optimization.
\newblock \emph{Manuscript submitted for publication to Information Sciences}, 2017.

\bibitem[Beume et~al.(2007)Beume, Naujoks, and Emmerich]{beume2007sms}
Nicola Beume, Boris Naujoks, and Michael Emmerich.
\newblock Sms-emoa: Multiobjective selection based on dominated hypervolume.
\newblock \emph{European Journal of Operational Research}, 181\penalty0
  (3):\penalty0 1653--1669, 2007.

\bibitem[Brockhoff et~al.(2015)Brockhoff, Tran, and
  Hansen]{brockhoff:hal-01146741}
Dimo Brockhoff, Thanh-Do Tran, and Nikolaus Hansen.
\newblock {Benchmarking Numerical Multiobjective Optimizers Revisited}.
\newblock In \emph{{Genetic and Evolutionary Computation Conference (GECCO
  2015)}}, Madrid, Spain, July 2015.
\newblock \doi{10.1145/2739480.2754777}.
\newblock URL \url{https://hal.inria.fr/hal-01146741}.

\bibitem[Cheng and Li(1999)]{6_p}
FY~Cheng and XS~Li.
\newblock Generalized center method for multiobjective engineering
  optimization.
\newblock \emph{Engineering Optimization}, 31\penalty0 (5):\penalty0 641--661,
  1999.

\bibitem[Coello et~al.(2002)Coello, Van~Veldhuizen, and
  Lamont]{coello2002evolutionary}
Carlos A~Coello Coello, David~A Van~Veldhuizen, and Gary~B Lamont.
\newblock \emph{Evolutionary algorithms for solving multi-objective problems},
  volume 242.
\newblock Springer, 2002.

\bibitem[Cust{\'o}dio et~al.(2011)Cust{\'o}dio, Madeira, Vaz, and
  Vicente]{custodio2011direct}
Ana~Lu{\'\i}sa Cust{\'o}dio, JF~Aguilar Madeira, A~Ismael~F Vaz, and
  Lu{\'\i}s~N Vicente.
\newblock Direct multisearch for multiobjective optimization.
\newblock \emph{SIAM Journal on Optimization}, 21\penalty0 (3):\penalty0
  1109--1140, 2011.

\bibitem[Deb(1999)]{15_p}
Kalyanmoy Deb.
\newblock Multi-objective genetic algorithms: Problem difficulties and
  construction of test problems.
\newblock \emph{Evolutionary computation}, 7\penalty0 (3):\penalty0 205--230,
  1999.

\bibitem[Deb et~al.(2002)Deb, Thiele, Laumanns, and Zitzler]{19_p}
Kalyanmoy Deb, Lothar Thiele, Marco Laumanns, and Eckart Zitzler.
\newblock Scalable multi-objective optimization test problems.
\newblock In \emph{Proceedings of the Congress on Evolutionary Computation
  (CEC-2002),(Honolulu, USA)}, pages 825--830. Proceedings of the Congress on
  Evolutionary Computation (CEC-2002),(Honolulu, USA), 2002.

\bibitem[Deb et~al.(2006)Deb, Sinha, and Kukkonen]{18_p}
Kalyanmoy Deb, Ankur Sinha, and Saku Kukkonen.
\newblock Multi-objective test problems, linkages, and evolutionary
  methodologies.
\newblock In \emph{Proceedings of the 8th annual conference on Genetic and
  evolutionary computation}, pages 1141--1148. ACM, 2006.

\bibitem[Fonseca and Fleming(1998)]{21_p}
Carlos~M Fonseca and Peter~J Fleming.
\newblock Multiobjective optimization and multiple constraint handling with
  evolutionary algorithms. i. a unified formulation.
\newblock \emph{Systems, Man and Cybernetics, Part A: Systems and Humans, IEEE
  Transactions on}, 28\penalty0 (1):\penalty0 26--37, 1998.

\bibitem[Hadka(2012)]{hadka2012moea}
D~Hadka.
\newblock Moea framework a free and open source java framework for
  multiobjective optimization, 2012.

\bibitem[Huband et~al.(2006)Huband, Hingston, Barone, and While]{25_p}
S.~Huband, P.~Hingston, L.~Barone, and L.~While.
\newblock A review of multiobjective test problems and a scalable test problem
  toolkit.
\newblock \emph{Evolutionary Computation, IEEE Transactions on}, 10\penalty0
  (5):\penalty0 477--506, Oct 2006.
\newblock ISSN 1089-778X.
\newblock \doi{10.1109/TEVC.2005.861417}.

\bibitem[Huband et~al.(2005)Huband, Barone, While, and Hingston]{24_p}
Simon Huband, Luigi Barone, Lyndon While, and Phil Hingston.
\newblock A scalable multi-objective test problem toolkit.
\newblock In \emph{Evolutionary multi-criterion optimization}, pages 280--295.
  Springer, 2005.

\bibitem[Hwang and Masud(1979)]{26_p}
C.-L. Hwang and A.~S.~MD. Masud.
\newblock Multiple objective decision making---methods and applications: A
  state-of-the-art survey.
\newblock \emph{Lecture Notes in Econom. Math. Systems}, 164, 1979.
\newblock Springer-Verlag, Berlin.

\bibitem[Knowles and Corne(2002)]{on_metric_performance}
J.~Knowles and D.~Corne.
\newblock On metrics for comparing nondominated sets.
\newblock In \emph{Evolutionary Computation, 2002. CEC '02. Proceedings of the
  2002 Congress on}, volume~1, pages 711--716, May 2002.
\newblock \doi{10.1109/CEC.2002.1007013}.

\bibitem[Knowles et~al.(2006)Knowles, Thiele, and
  Zitzler]{knowles_tutorial_indicator}
J.D. Knowles, L.~Thiele, and E.~Zitzler.
\newblock A tutorial on the performance assessment of stochastic
  multi-objective optimizers.
\newblock TIK-Report 214, Computer Engineering and Networks Laboratory, ETH
  Zurich, Gloriastrasse 35, ETH-Zentrum, 8092 Zurich, Switzerland, February
  2006.

\bibitem[Kolda et~al.(2003)Kolda, Lewis, and Torczon]{31_p}
Tamara~G Kolda, Robert~Michael Lewis, and Virginia Torczon.
\newblock Optimization by direct search: New perspectives on some classical and
  modern methods.
\newblock \emph{SIAM review}, 45\penalty0 (3):\penalty0 385--482, 2003.

\bibitem[Liuzzi et~al.(2003)Liuzzi, Lucidi, Parasiliti, and Villani]{33_p}
Giampaolo Liuzzi, Stefano Lucidi, Francesco Parasiliti, and Marco Villani.
\newblock Multiobjective optimization techniques for the design of induction
  motors.
\newblock \emph{IEEE Transactions on Magnetics}, 39\penalty0 (3):\penalty0
  1261--1264, 2003.

\bibitem[Mor{\'e} and Wild(2009)]{more2009benchmarking}
Jorge~J Mor{\'e} and Stefan~M Wild.
\newblock Benchmarking derivative-free optimization algorithms.
\newblock \emph{SIAM Journal on Optimization}, 20\penalty0 (1):\penalty0
  172--191, 2009.

\bibitem[Nocedal and Wright(2006)]{39_p}
J.~Nocedal and Stephen~J Wright.
\newblock \emph{Numerical optimization}, volume~2.
\newblock Springer-Verlag, Berlin, 2 edition, 2006.

\bibitem[Vicente and Cust{\'o}dio(2012)]{48_p}
Lu{\'\i}s~N Vicente and AL~Cust{\'o}dio.
\newblock Analysis of direct searches for discontinuous functions.
\newblock \emph{Mathematical programming}, 133\penalty0 (1-2):\penalty0
  299--325, 2012.

\bibitem[Zitzler et~al.(2000)Zitzler, Deb, and Thiele]{49_p}
Eckart Zitzler, Kalyanmoy Deb, and Lothar Thiele.
\newblock Comparison of multiobjective evolutionary algorithms: Empirical
  results.
\newblock \emph{Evolutionary computation}, 8\penalty0 (2):\penalty0 173--195,
  2000.

\bibitem[Zitzler et~al.(2003)Zitzler, Thiele, Laumanns, Fonseca, and
  Da~Fonseca]{zitzler2003performance}
Eckart Zitzler, Lothar Thiele, Marco Laumanns, Carlos~M Fonseca, and
  Viviane~Grunert Da~Fonseca.
\newblock Performance assessment of multiobjective optimizers: an analysis and
  review.
\newblock \emph{IEEE Transactions on Evolutionary Computation}, 7\penalty0
  (2):\penalty0 117--132, 2003.

\end{thebibliography}
\end{document}